# EQUILIBRIUM FOR FRAGMENTATION WITH IMMIGRATION

By Bénédicte Haas

*Université Pierre et Marie Curie et C.N.R.S. UMR 7599*

This paper introduces stochastic processes that describe the evolution of systems of particles in which particles immigrate according to a Poisson measure and split according to a self-similar fragmentation. Criteria for existence and absence of stationary distributions are established and uniqueness is proved. Also, convergence rates to the stationary distribution are given. Linear equations which are the deterministic counterparts of fragmentation with immigration processes are next considered. As in the stochastic case, existence and uniqueness of solutions, as well as existence and uniqueness of stationary solutions, are investigated.

**1. Introduction.** The aim of this paper is to study random and deterministic models that describe the evolution of systems of particles in which two independent phenomena take place: immigration and fragmentation of particles. Particles immigrate and split into smaller particles, which, in turn, continue splitting, at rates that depend on their mass. Such a situation occurs, for example, in grinding lines [1, 23] where macroscopic blocks are continuously placed in tumbling ball mills that reduce them to microscopic fragments. These microscopic fragments then undergo a chemical process to extract the minerals. In such systems, one may expect to attain an equilibrium, as the immigration may compensate for the fragmentation of particles. The investigation of existence and uniqueness of such stationary state, as well as convergence to the stationary state, is one of the main points of interest of this paper. It will be undertaken both in random and deterministic settings.

We first introduce continuous-time *fragmentation with immigration Markov processes*. Roughly, their dynamics are described as follows. The immigration is coded by a Poisson measure with intensity $I(d\mathbf{s})\,dt$, $t \geq 0$, where $I$ is a measure supported on $\mathcal{D}$, the set of decreasing sequences $\mathbf{s} = (s_j, j \geq 1)$ that converge to 0. That is, if $(\mathbf{s}(t_i), t_i)$ denotes the atoms of this Poisson measure, a group of particles with masses $(s_1(t_i), s_2(t_i), \dots)$ immigrates at time $t_i$ for each $t_i \geq 0$. We further impose that $I$ integrates $\sum_{j \geq 1}(s_j \wedge 1)$, which means that the total mass of immigrants on a finite time interval is finite a.s. The particles fragment independently of the immigration, according to a "self-similar fragmentation with index $\alpha \in \mathbb{R}$" as introduced by Bertoin in [6, 7]. This means that each particle splits independently of others with a rate proportional to its mass to the power $\alpha$ and that the resulting particles continue splitting with the same rules. Rigorous definitions are given in Sections 1.1 and 1.2 below. Some examples of such processes arise from classical stochastic processes, as a Brownian motion with positive drift. This is detailed in Section 4.

Let $FI$ denote a fragmentation with immigration process. Our first purpose is to know whether it is possible to find a *stationary distribution* for $FI$. Let us mention here that, until now, an equilibrium could only be obtained for fragmentation with *coagulation* processes. See, for example, [4, 14, 16].

Under some conditions that depend both on the dynamics of the fragmentation and on the immigration, we construct a random variable $\mathbf{U}_{\text{stat}}$ in $\mathcal{D}$ whose distribution is stationary for $FI$. Let $\alpha_I$ be the $I$-dependent parameter defined by

$$(1) \qquad \alpha_I := -\sup\left\{a > 0 : \int_{\mathcal{D}} s_1^a \mathbf{1}_{\{s_1 \geq 1\}} I(d\mathbf{s}) < \infty\right\}.$$

When $\alpha_I < 0$, we obtain that the stationary state $\mathbf{U}_{\text{stat}}$ exists as soon as the index of self-similarity $\alpha$ is larger than $\alpha_I$ and that there is no stationary distribution when $\alpha$ is smaller than $\alpha_I$. In this latter case, the particles with mass larger than 1, which split slower when $\alpha$ is smaller, do not split fast enough to compensate the immigration of large particles, which therefore accumulate. In others words, too many large particles are brought in the ball mill which is not able to grind them fast enough. These results are made precise in Theorems 7, 8 and 9, Section 2, where we also study whether $\mathbf{U}_{\text{stat}}$ is in $l^p$, $p \geq 0$. In addition, the stationary solution is proved unique.

It is easily checked from the construction of $\mathbf{U}_{\text{stat}}$ that

$$FI(t) \overset{\text{law}}{\to} \mathbf{U}_{\text{stat}}$$

as soon as the stationary distribution exists and that this convergence holds independently of the initial distribution. One standard problem is to investigate the rate of convergence to this stationary state. Our approach is based on a coupling method. This provides rates of convergence that differ significantly according as $\alpha < 0$, $\alpha = 0$ or $\alpha > 0$: one obtains that the convergence takes place at a geometric rate when $\alpha = 0$, at rate $t^{-1/\alpha}$ when $\alpha > 0$, whereas the rate of convergence depends both on $I$ and $\alpha$ when $\alpha < 0$.

We next turn to deterministic models, namely, *fragmentation with immigration equations*. Roughly, these equations are obtained by adding an immigration term to a family of well-known fragmentation equations with mass loss [17, 18, 24]: we consider that particles with mass in the interval $(x, x+dx)$ arrive at rate $\mu_I(dx)$ which is defined from $I$ by

$$\int_0^\infty f(x)\mu_I(dx) := \int_\mathcal{D} \sum_{j\geq 1} f(s_j) I(d\mathbf{s}),$$

for all positive measurable functions $f$. Solutions to the fragmentation with immigration equation do not always exist. We give conditions for existence and then show uniqueness. The obtained solution is closely related to the stochastic model $(FI(t), t \geq 0)$: it is—in a sense to be specified—related to the expectations of the random measures $\sum_{k\geq 1} \delta_{FI_k(t)}$, $t \geq 0$. In this deterministic setting, one may also expect the existence of stationary solutions. Provided the average mass immigrated by unit time is finite, we construct explicitly a stationary solution which is proved unique. Note that here the hypothesis for existence only involves $I$, not $\alpha$, contrary to the stochastic case.

This paper is organized as follows. In the remainder of this section we first review the definition and some properties of self-similar fragmentations (Section 1.1), then we set down the definition of fragmentation with immigration processes (Section 1.2). The study of existence and uniqueness of a stationary distribution is undertaken in Section 2, where we also give criteria for existence of a stationary distribution for more general Markov processes with immigration. In Section 3 we investigate the rate of convergence to the stationary distribution. Section 4 is devoted to examples of fragmentation with immigration processes constructed from Brownian motions with positive drift. Section 5 concerns the fragmentation with immigration equation.

### 1.1. Self-similar fragmentations.

*State space.* We endow the state space

$$\mathcal{D} = \left\{ \mathbf{s} = (s_j)_{j\geq 1} : s_1 \geq s_2 \geq \cdots \geq 0, \lim_{j\to\infty} s_j = 0 \right\}$$

with the uniform distance

$$d(\mathbf{s}, \mathbf{s}') := \sup_{j\geq 1} |s_j - s_j'|.$$

Clearly, as $n \to \infty$, $d(\mathbf{s}, \mathbf{s}^n) \to 0$ is equivalent to $s_j^n \to s_j$ for all $j \geq 1$, which, in turn, is equivalent to $\sum_{j\geq 1} f(s_j^n) \to \sum_{j\geq 1} f(s_j)$ for all continuous functions $f$ with compact support in $(0, \infty)$. Hence, $\mathcal{D}$ identifies with the set of Radon counting measures on $(0, \infty)$ with bounded support endowed with the topology of vague convergence through the homeomorphism

$$\mathbf{s} \in \mathcal{D} \mapsto \sum_{j \geq 1} \delta_{s_j} \mathbf{1}_{\{s_j > 0\}}.$$

With a slight abuse of notation, we also call $\mathbf{s}$ the measure $\sum_{j \geq 1} \delta_{s_j} \mathbf{1}_{\{s_j > 0\}}$. It is then natural to denote by "$\mathbf{s} + \mathbf{s}'$" the decreasing rearrangement of the concatenation of sequences $\mathbf{s}, \mathbf{s}'$ and by $\langle \mathbf{s}, f \rangle$ the sum $\sum_{j \geq 1} f(s_j) \mathbf{1}_{\{s_j > 0\}}$. More generally, we denote by "$\sum_{i \geq 1} \mathbf{s}^i$" the measure $\sum_{i \geq 1} \sum_{j \geq 1} \delta_{s_j^i} \mathbf{1}_{\{s_j^i > 0\}}$. This point measure does not necessarily correspond to a sequence in $\mathcal{D}$, but when it does, it represents the decreasing rearrangement of the concatenation of sequences $\mathbf{s}^1, \mathbf{s}^2, \ldots$.

For all $p \geq 0$, let $l^p$ be the subset of $\mathcal{D}$ of sequences $s_1 \geq s_2 \geq \cdots \geq 0$ such that $\sum_{j \geq 1} s_j^p < \infty$, endowed with the topology of $\mathcal{D}$. When $p = 0$, we use the convention $0^0 = 0$, which means that $l^0$ is the space of sequences with at most a finite number of nonzero terms. Let also $\mathcal{D}_1$ be the subset of $\mathcal{D}$ of sequences such that $\sum_{j \geq 1} s_j \leq 1$. Clearly, $l^p \subset l^{p'}$ when $p \leq p'$ and $\mathcal{D}_1 \subset l^1$. At last, set $\mathbf{0} := (0, 0, \ldots)$.

*Self-similar fragmentations.*

DEFINITION 1. A standard self-similar fragmentation $(F(t), t \geq 0)$ with index $\alpha \in \mathbb{R}$ is a $\mathcal{D}_1$-valued Markov process continuous in probability such that:

- $F(0) = (1, 0, \ldots)$,
- for each $t_0 \geq 0$, conditionally on $F(t_0) = (s_1, s_2, \ldots)$, the process $(F(t + t_0), t \geq 0)$ has the same law as the process obtained for each $t \geq 0$ by ranking in the decreasing order the components of sequences $s_1 F^{(1)}(s_1^\alpha t)$, $s_2 F^{(2)}(s_2^\alpha t), \ldots$, where the $F^{(j)}$'s are independent copies of $F$.

This means that the particles present at a time $t_0$ evolve independently and that the evolution process of a particle with mass $m$ has the same distribution as $m$ times the process starting from a particle with mass 1, up to the time change $t \mapsto t m^\alpha$. According to [3] and [7], a self-similar fragmentation is Feller—hence, possesses a càdlàg version which we shall always consider—and its distribution is characterized by a 3-tuple $(\alpha, c, \nu)$: $\alpha$ is the index of self-similarity, $c \geq 0$ an *erosion* coefficient and $\nu$ a *dislocation measure*, which is a sigma-finite nonnegative measure on $\mathcal{D}_1$ that does not charge $(1, 0, \ldots)$ and satisfies

$$\int_{\mathcal{D}_1} (1 - s_1) \nu(d\mathbf{s}) < \infty.$$

Roughly speaking, the erosion is a deterministic continuous phenomenon and the dislocation measure describes the rates of sudden dislocations: a fragment with mass $m$ splits into fragments with masses $m\mathbf{s}, \mathbf{s} \in \mathcal{D}_1$, at rate $m^\alpha \nu(d\mathbf{s})$. In case $\nu(\mathcal{D}_1) < \infty$ and $c = 0$, this means that a particle with mass $m$ splits after a time $T$ with an exponential law with parameter $m^\alpha \nu(\mathcal{D}_1)$ into particles with masses $m\mathbf{s}$, where $\mathbf{s}$ is distributed according to $\nu(\cdot)/\nu(\mathcal{D}_1)$ and is independent of $T$. For more details on these fundamental properties of self-similar fragmentations, we refer to [3, 6, 7].

DEFINITION 2. For any random $\mathbf{u} \in \mathcal{D}$, a fragmentation process $(\alpha, c, \nu)$, starting from $\mathbf{u}$, is defined by

$$(2) \qquad F^{(\mathbf{u})}(t) := \sum_{j \geq 1} (u_j F^{(j)}(u_j^\alpha t)), \qquad t \geq 0,$$

where the $F^{(j)}$'s are i.i.d. copies of a standard $(\alpha, c, \nu)$-fragmentation $F$, independent of $\mathbf{u}$.

Clearly, $F^{(\mathbf{u})}(t) \in \mathcal{D}$ for all $t \geq 0$ and, according to the branching property of $F$, $F^{(\mathbf{u})}$ is Markov. It is plain that such a fragmentation process converges a.s. to $\mathbf{0}$ as $t \to \infty$, provided $\nu(\mathcal{D}_1) \neq 0$. We shall denote in the sequel by $F_1^{(\mathbf{u})}(t) \geq \cdots \geq F_k^{(\mathbf{u})}(t) \geq \cdots$ the components of the sequence $F^{(\mathbf{u})}(t)$.

We now review some facts about standard $(\alpha, c, \nu)$-fragmentations that we will need. In the remainder of this section, $F$ denotes a standard $(\alpha, c, \nu)$-fragmentation.

*Tagged particle.* We are interested in the evolution process $\lambda$ of the mass of a particle tagged at random in the fragmentation. To construct this process, we recall that one may always suppose that $F$ is built from some family $(G(t), t \geq 0)$ of nested open sets of $(0, 1)$ so that $F(t)$ is the ordered sequence of lengths of the interval components of $G(t)$, $t \geq 0$ (see [3, 7]). Let then $U$ be uniformly distributed on $(0, 1)$, independent of $G$, and call $\lambda(t)$ the length of the interval component of $G(t)$ containing $U$. When such interval does not exist, set $\lambda(t) := 0$. The main point of interest of such approach is that the distribution of $\lambda$ is well known.

First, when $\alpha = 0$, Bertoin [6] shows that $\lambda \stackrel{\text{law}}{=} \exp(-\xi(\cdot))$, where $\xi$ is a subordinator (i.e., a right-continuous increasing process with values in $[0, \infty]$ and with stationary and independent increments on the interval $\{t : \xi(t) < \infty\}$), with Laplace exponent $\phi$ given by

$$(3) \qquad \phi(q) := c(q+1) + \int_{\mathcal{D}_1} \left(1 - \sum_{j \geq 1} s_j^{1+q}\right) \nu(d\mathbf{s}), \qquad q \geq 0.$$

We recall that $\phi$ characterizes $\xi$, since $E[\exp(-q\xi(t))] = \exp(-t\phi(q))$ for all $t, q \geq 0$ (for background on subordinators, we refer to [5], Chapter III). When $c > 0$ or $\nu(\sum_{j \geq 1} s_j < 1) > 0$, one sees that the subordinator $\xi$ is killed at rate $k = \phi(0) > 0$: that is, there exists a subordinator $\overline{\xi}$ with Laplace exponent $\overline{\phi} = \phi - k$ and an exponential r.v. $\mathbf{e}(k)$ with parameter $k$, independent of $\overline{\xi}$, such that

$$\xi(t) = \overline{\xi}(t)\mathbf{1}_{\{t < \mathbf{e}(k)\}} + \infty \mathbf{1}_{\{t \geq \mathbf{e}(k)\}}$$

for all $t \geq 0$.

When $\alpha \in \mathbb{R}$, Bertoin [7] shows that $\lambda \stackrel{\text{law}}{=} \exp(-\xi(\rho(\cdot)))$, where $\xi$ is the same subordinator as above and $\rho$ is the time-change

(4) $$\rho(t) := \inf\left\{u \geq 0 : \int_0^u \exp(\alpha \xi(r))\, dr > t\right\}, \qquad t \geq 0.$$

On the other hand, the construction of $\lambda$ implies that, conditionally on $F$, $\lambda(t) = F_k(t)$ with probability $F_k(t)$, $k \geq 1$, and that $\lambda(t) = 0$ with probability $1 - \sum_{k \geq 1} F_k(t)$. Hence,

(5) $$\sum_{k \geq 1} f(F_k(t)) = E[f(\exp(-\xi(\rho(t))))\exp(\xi(\rho(t)))|F]$$

for every positive measurable function $f$ supported on a compact subset of $(0, \infty)$ (with the convention $0 \times \infty = 0$), and, in particular,

(6) $$E\left[\sum_{k \geq 1} f(F_k(t))\right] = E[f(\exp(-\xi(\rho(t))))\exp(\xi(\rho(t)))].$$

*Formation of dust when $\alpha < 0$.* When the index of self-similarity $\alpha$ is negative, for all dislocation measures $\nu$, the total mass $\sum_{k \geq 1} F_k(t)$ of the fragmentation $F$ decreases as time passes to reach 0 in finite time even if there is no erosion ($c = 0$) and no mass is lost within sudden dislocations ($\nu(\sum_{j \geq 1} s_j < 1) = 0$). This is due to an intensive fragmentation of small particles which reduces macroscopic particles to an infinite number of zero-mass particles or *dust*. To say this precisely, introduce

(7) $$\tau := \inf\left\{t \geq 0 : \sum_{k \geq 1} F_k(t) = 0\right\}$$

the first time at which the total mass reaches 0. According to Proposition 14 in [18], there exist $C, C'$ some positive finite constants such that, for any $t \geq 0$,

(8) $$P(\tau > t) \leq C \exp(-C' t^\Gamma),$$

where $\Gamma$ is a $(c, \nu)$-dependent parameter defined by

(9) $$\Gamma := \begin{cases} (1-\lambda)^{-1}, & \text{when } \phi(q) - cq \\ & \text{varies regularly with index } 0 < \lambda < 1 \text{ as } q \to \infty, \\ 1, & \text{otherwise.} \end{cases}$$

Note that $E[\tau]<\infty$. This phenomenon of formation of dust does not occur when $\alpha \geq 0$: if no mass is lost by erosion or within sudden dislocations, then $\sum_{k\geq 1} F_k(t) = 1$ a.s. for all $t \geq 0$.

1.2. *Fragmentation with immigration processes.* As said previously, the immigration and fragmentation phenomena occur independently. The immigration is coded by a Poisson measure on $l^1 \times [0,\infty)$ with an intensity $I(d\mathbf{s})\,dt$ such that

$$\text{(H1)} \qquad \int_{l^1} \sum_{j\geq 1}(s_j \wedge 1) I(d\mathbf{s}) < \infty$$

and we call such measure $I$ an *immigration measure*. The hypothesis (H1) implies that the total mass of particles that have immigrated during a time $t$ is almost surely finite (for an introduction to Poisson measures, we refer to [21]). On the other hand, the particles fragment according to a self-similar fragmentation $(\alpha, c, \nu)$.

DEFINITION 3. Let $\mathbf{u}$ be a random sequence of $\mathcal{D}$ and let $((\mathbf{s}(t_i), t_i), i \geq 1)$ be the atoms of a Poisson measure with intensity $I(d\mathbf{s})\,dt$ independent of $\mathbf{u}$. Then, conditionally on $\mathbf{u}$ and $((\mathbf{s}(t_i), t_i), i \geq 1)$, let $F^{(\mathbf{u})}, F^{(\mathbf{s}(t_i))}, i \geq 1$, be independent fragmentation processes $(\alpha, c, \nu)$ starting, respectively, from $\mathbf{u}, \mathbf{s}(t_1), \mathbf{s}(t_2), \ldots$. With probability one, the sum

$$FI^{(\mathbf{u})}(t) := F^{(\mathbf{u})}(t) + \sum_{t_i \leq t} F^{(\mathbf{s}(t_i))}(t - t_i)$$

belongs to $\mathcal{D}$ for all $t \geq 0$, and the process $FI^{(\mathbf{u})}$ is called a fragmentation with immigration process with parameters $(\alpha, c, \nu, I)$ starting from $\mathbf{u}$.

One may be troubled by conditioning on the value of $((\mathbf{s}(t_i), t_i), i \geq 1)$, as it may have 0 probability. If so, note that the family $F^{(\mathbf{s}(t_i))}, i \geq 1$, is actually constructed from the Poisson measure $((\mathbf{s}(t_i), t_i), i \geq 1)$ and an independent family $F^{(i,j)}, i, j \geq 1$, of i.i.d. standard $(\alpha, c, \nu)$-fragmentations, through the formula $F^{(\mathbf{s}(t_i))}(t) = \sum_{i,j\geq 1} s_j(t_i) F^{(i,j)}((s_j(t_i))^{\alpha} t)$.

The reason why $\sum_{t_i \leq t} F^{(\mathbf{s}(t_i))}(t - t_i) \in \mathcal{D}$ a.s. is that $\sum_{t_i \leq t} \sum_{j \geq 1} s_j(t_i) < \infty$ [by hypothesis (H1)] and then that $\sum_{t_i \leq t} F^{(\mathbf{s}(t_i))}(t - t_i) \in l^1$, since $\sum_{k\geq 1} F_k^{(\mathbf{s}(t_i))}(t - t_i) \leq \sum_{j\geq 1} s_j(t_i)$. Note also that when $p \geq 1$, $FI^{(\mathbf{u})} \in l^p$ as soon as $\mathbf{u} \in l^p$.

In this definition, the sequence $\mathbf{u}$ represents the masses of particles present at time 0 and at each time $t_i \geq 0$, some particles of masses $\mathbf{s}(t_i)$ immigrate. At time $t$, two families of particles are then present: those resulting from the fragmentation of $\mathbf{u}$ during a time $t$ and those resulting from the fragmentation of $\mathbf{s}(t_i)$ during a time $t - t_i$, $t_i \leq t$.

It is easy to see that the process $FI^{(\mathbf{u})}$ is Markov with the Feller property (cf. the proof of Proposition 1.1 in [3]). Hence, we may and will always consider càdlàg versions of $FI^{(\mathbf{u})}$.

In the rest of this paper, we denote by $FI$ a fragmentation with immigration $(\alpha, c, \nu, I)$ (without any specified starting point) and we always exclude the trivial cases $\nu = 0$ or $I = 0$.

REMARK. One may wonder why we do not more generally consider some fragmentation with immigration processes with values in $\mathcal{R}$, the set of Radon point measures on $(0, \infty)$. Indeed, for all (random) $\mathbf{u} \in \mathcal{R}$ and all $t \geq 0$, it is always possible to define the point measure

$$(10) \qquad FI^{(\mathbf{u})}(t) := F^{(\mathbf{u})}(t) + \sum_{t_i \leq t} F^{(\mathbf{s}(t_i))}(t - t_i), \qquad t \geq 0,$$

where $F^{(\mathbf{u})}(t)$ is defined similarly to (2) and is independent of $F^{(\mathbf{s}(t_i))}, i \geq 1$, some independent fragmentations $(\alpha, c, \nu)$ starting, respectively, from $\mathbf{s}(t_1)$, $\mathbf{s}(t_2), \ldots$. The sum involving the terms $F^{(\mathbf{s}(t_i))}(t - t_i)$, $t_i \leq t$, is in $\mathcal{D}$, as noticed in the Definition 3 above. The issue is that, in general, starting from some $\mathbf{u} \in \mathcal{R} \backslash \mathcal{D}$, the measures $F^{(\mathbf{u})}(t)$ do not necessarily belong to $\mathcal{R}$, as the masses of the initial particles may accumulate in some bounded interval $(a, b)$ after fragmentation.

As an example, starts from $\mathbf{u} = \sum_{i \geq 1} \delta_i$ and fix $\alpha > 0$. For each $i$, tag a particle at random in the fragmentation issued from the particle with mass $i$, as explained in the previous section. At time $t$, this tagged particle is distributed as $i \exp(-\xi^{(i)}(\rho^{(i)}(i^\alpha t)))$, where the $\xi^{(i)}$'s are i.i.d. subordinators with Laplace exponent (3) and $\rho^{(i)}$ the corresponding time changes (4). According to the Borel–Cantelli lemma, the number of tagged particles belonging to some interval $(a, b)$ at time $t$ is then a.s. infinite [and, therefore, $F^{(\mathbf{u})}(t) \notin \mathcal{R}$] as soon as $\sum_{i \geq 1} P(a < i \exp(-\xi(\rho(i^\alpha t))) < b) = \infty$. In [9], Bertoin and Caballero show that for most of subordinators (and, therefore, for most of dislocation measures) $i \exp(-\xi(\rho(i^\alpha t)))$ has a nontrivial limiting distribution as $i \to \infty$ when $\alpha > 0$. In such cases, the above sum of probabilities is infinite for some well-chosen intervals $(a, b)$ and then $F^{(\mathbf{u})}(t) \notin \mathcal{R}$.

That is why we study fragmentation with immigration processes on $\mathcal{D}$. However, in Section 5, we shall use some of these measures $FI^{(\mathbf{u})}(t)$, $\mathbf{u} \in \mathcal{R}$, and we give (Proposition 15) some sufficient conditions on $\mathbf{u}$ and $\alpha$ for $F^{(\mathbf{u})}(t)$ [equivalently, $FI^{(\mathbf{u})}(t)$] to be a.s. Radon at fixed time $t$. These conditions do not ensure that the process $FI^{(\mathbf{u})}$ is $\mathcal{R}$-valued, as we do not know if a.s., for *all* $t$, $FI^{(\mathbf{u})}(t) \in \mathcal{R}$.

**2. Existence and uniqueness of the stationary distribution.** This section is devoted to the existence and uniqueness of a stationary distribution for *FI* and to properties of the stationary state, when it exists. We begin by establishing some criteria for existence and uniqueness of a stationary distribution, which are available for a class of Markov processes with immigration including fragmentation with immigration processes. This is undertaken in Section 2.1 where we more specifically obtain an explicit construction of a stationary state. We then apply these results to fragmentation with immigration processes (Section 2.2).

From now on, for any r.v. $X$, $\mathcal{L}(X)$ denotes the distribution of $X$.

2.1. *The candidate for a stationary distribution for Markov processes with immigration.* Recall that $\mathcal{R}$ denotes the set of Radon point measures on $(0, \infty)$ and equip it with the topology of vague convergence. We first consider $\mathcal{R}$-valued Markov processes with some superposition property and then extend the results to a larger class of Markov processes.

Let $X$ be an $\mathcal{R}$-valued Markov process that satisfies the following *superposition property*: for all $\mathbf{u}, \mathbf{v} \in \mathcal{R}$, the sum of two independent processes $X^{(\mathbf{u})}$ and $X^{(\mathbf{v})}$ starting, respectively, from $\mathbf{u}$ and $\mathbf{v}$ is distributed as $X^{(\mathbf{u}+\mathbf{v})}$. A moment of thought shows that this is equivalent to $\sum_{i \geq 1} X^{(\mathbf{u}^i)} \stackrel{\text{law}}{=} X^{(\sum_{i \geq 1} \mathbf{u}^i)}$ for all sequences $(\mathbf{u}^i, i \geq 1)$ such that $\sum_{i \geq 1} \mathbf{u}^i \in \mathcal{R}$ a.s., where $X^{(\mathbf{u}^1)}, X^{(\mathbf{u}^2)}, \ldots$ are independent processes, starting, respectively, from $\mathbf{u}^1, \mathbf{u}^2, \ldots$. Consider then $I$, a nonnegative $\sigma$-finite measure on $\mathcal{R}$, and let $((\mathbf{s}(t_i), t_i), i \geq 1)$ be the atoms of a Poisson measure with intensity $I(d\mathbf{s})\, dt, t \geq 0$. Conditionally on this Poisson measure, let $X^{(\mathbf{s}(t_i))}$ be independent versions of $X$, starting, respectively, from $\mathbf{s}(t_1), \mathbf{s}(t_2), \ldots$. In order to define an $X$-*process with immigration*, we need and will suppose in this section that a.s.

$$\sum_{t_i \leq t} X^{(\mathbf{s}(t_i))}(t - t_i) \in \mathcal{R} \qquad \text{for all } t \geq 0.$$

In particular, this holds when $X$ is a fragmentation process and $I$ an immigration measure, as explained just after Definition 3. More generally, still supposing that $X$ is a fragmentation, one easily checked that it holds as soon as $I$ integrates $\mathbf{1}_{\{s_1 > \varepsilon\}}$ for all $\varepsilon > 0$, that is, as soon as the number of particles of mass larger than $\varepsilon$ immigrating in finite time is finite.

DEFINITION 4. For every random $\mathbf{u} \in \mathcal{R}$, let $X^{(\mathbf{u})}$ be a version of $X$ starting from $\mathbf{u}$ and consider $((X^{(\mathbf{r}(v_i))}, v_i), i \geq 1)$ a version of $((X^{(\mathbf{s}(t_i))}, t_i), i \geq 1)$ independent of $X^{(\mathbf{u})}$. Then, the process defined by

(11) $$XI^{(\mathbf{u})}(t) := X^{(\mathbf{u})}(t) + \sum_{v_i \leq t} X^{(\mathbf{r}(v_i))}(t - v_i), \qquad t \geq 0,$$

is an $\mathcal{R}$-valued Markov process and is called $X$-*process with immigration starting from* $\mathbf{u}$.

We point out that the Markov property of $XI$ results both from the Markov property and from the superposition property of $X$. A moment of reflection shows that the law of the point measure

$$\mathbf{U}_{\text{stat}} := \sum_{t_i \geq 0} X^{(\mathbf{s}(t_i))}(t_i) \tag{12}$$

is a natural candidate for a stationary distribution for $XI$ [in some sense, it is the limit as $t \to \infty$ of $XI^{(\mathbf{0})}(t)$], provided that it belongs to $\mathcal{R}$. The problem is that it does not necessarily belong to $\mathcal{R}$, as the components of $\mathbf{U}_{\text{stat}}$ may accumulate in some bounded interval $(a, b)$.

LEMMA 5. (i) *If $\mathbf{U}_{\text{stat}} \in \mathcal{R}$ a.s., then the distribution $\mathcal{L}(\mathbf{U}_{\text{stat}})$ is a stationary distribution for $XI$ and for any random $\mathbf{u} \in \mathcal{R}$ such that $X^{(\mathbf{u})}(t) \overset{\text{P}}{\to} \mathbf{0}$ as $t \to \infty$,*

$$XI^{(\mathbf{u})}(t) \overset{\text{law}}{\to} \mathbf{U}_{\text{stat}} \qquad \text{as } t \to \infty.$$

(ii) *If $P(\mathbf{U}_{\text{stat}} \notin \mathcal{R}) > 0$, then there exists no stationary distribution for $XI$ and if $P(\mathbf{U}_{\text{stat}} \notin \mathcal{D}) > 0$, then there exists no stationary distribution on $\mathcal{D}$ for $XI$.*

PROOF. (i) Assume $\mathbf{U}_{\text{stat}} \in \mathcal{R}$ a.s. and consider a version $XI^{(\mathbf{U}_{\text{stat}})}$ of the $X$-process with immigration starting from $\mathbf{U}_{\text{stat}}$. We want to prove that $XI^{(\mathbf{U}_{\text{stat}})}(t) \overset{\text{law}}{=} \mathbf{U}_{\text{stat}}$ for every $t \geq 0$. So fix $t > 0$. By definition of $XI$ and using the Markov and superposition properties of $X$, we see that there exists $((X^{(\mathbf{r}(v_i))}, v_i), i \geq 1)$ an independent copy of $((X^{(\mathbf{s}(t_i))}, t_i), i \geq 1)$ such that

$$XI^{(\mathbf{U}_{\text{stat}})}(t) \overset{\text{law}}{=} \sum_{t_i \geq 0} X^{(\mathbf{s}(t_i))}(t_i + t) + \sum_{v_i \leq t} X^{(\mathbf{r}(v_i))}(t - v_i).$$

By independence of $((\mathbf{r}(v_i), v_i), i \geq 1)$ and $((\mathbf{s}(t_i), t_i), i \geq 1)$, the concatenation of

$$((\mathbf{r}(v_i), t - v_i), v_i \leq t) \quad \text{and} \quad ((\mathbf{s}(t_i), t_i + t), i \geq 1)$$

has the same law as $((\mathbf{s}(t_i), t_i), i \geq 1)$. Hence,

$$XI^{(\mathbf{U}_{\text{stat}})}(t) \overset{\text{law}}{=} \sum_{t_i \geq 0} X^{(\mathbf{s}(t_i))}(t_i) = \mathbf{U}_{\text{stat}}.$$

Similarly, one obtains that, for all $t \geq 0$,

$$XI^{(\mathbf{u})}(t) \overset{\text{law}}{=} X^{(\mathbf{u})}(t) + \sum_{v_i \leq t} X^{(\mathbf{r}(v_i))}(v_i), \tag{13}$$

where $((X^{(\mathbf{r}(v_i))}, v_i), i \geq 1)$ is distributed as $((X^{(\mathbf{s}(t_i))}, t_i), i \geq 1)$ and is independent of $X^{(\mathbf{u})}$. Suppose now that $X^{(\mathbf{u})}(t) \xrightarrow{\mathrm{P}} \mathbf{0}$ as $t \to \infty$. Clearly,

$$\sum_{v_i \leq t} X^{(\mathbf{r}(v_i))}(v_i) \xrightarrow[t \to \infty]{\mathrm{a.s.}} \sum_{v_i \geq 0} X^{(\mathbf{r}(v_i))}(v_i)$$

and, therefore,

$$X^{(\mathbf{u})}(t) + \sum_{v_i \leq t} X^{(\mathbf{r}(v_i))}(v_i) \xrightarrow{\mathrm{P}} \sum_{v_i \geq 0} X^{(\mathbf{r}(v_i))}(v_i) \qquad \text{as } t \to \infty.$$

Since the limit here is distributed as $\mathbf{U}_{\mathrm{stat}}$ and since (13) holds, one has $XI^{(\mathbf{u})}(t) \xrightarrow{\mathrm{law}} \mathbf{U}_{\mathrm{stat}}$.

(ii) Suppose that there exists a stationary distribution $\mathcal{L}_{\mathrm{stat}}$. Our aim is to show that $P(\mathbf{U}_{\mathrm{stat}} \notin \mathcal{R}) = 0$. To do so, let $XI^{(\mathcal{L}_{\mathrm{stat}})}$ be an $X$-process with immigration starting from an initial sequence distributed according to $\mathcal{L}_{\mathrm{stat}}$. Replacing $\mathbf{u}$ by $XI^{(\mathcal{L}_{\mathrm{stat}})}(0)$ in (13), we get

$$XI^{(\mathcal{L}_{\mathrm{stat}})}(0) \stackrel{\mathrm{law}}{=} X^{(XI^{(\mathcal{L}_{\mathrm{stat}})}(0))}(t) + \sum_{t_i \leq t} X^{(\mathbf{s}(t_i))}(t_i).$$

Introduce then, for any $0 < a < b < \infty$, the event

$$E_{a,b} := \left\{ \sum_{t_i \geq 0} \langle X^{(\mathbf{s}(t_i))}(t_i), \mathbf{1}_{(a,b)} \rangle = \infty \right\}$$

and fix some $N > 0$. The identity in law obtained above yields

$$P(\langle XI^{(\mathcal{L}_{\mathrm{stat}})}(0), \mathbf{1}_{(a,b)} \rangle < N)$$
$$\leq P\left( \sum_{t_i \leq t} \langle X^{(\mathbf{s}(t_i))}(t_i), \mathbf{1}_{(a,b)} \rangle < N \right)$$
$$\leq P\left( \sum_{t_i \leq t} \langle X^{(\mathbf{s}(t_i))}(t_i), \mathbf{1}_{(a,b)} \rangle < N, E_{a,b} \right) + P(\Omega \backslash E_{a,b}).$$

The first probability in this latter sum converges to 0 as $t \to \infty$ by definition of $E_{a,b}$ and, therefore,

$$P(\langle XI^{(\mathcal{L}_{\mathrm{stat}})}(0), \mathbf{1}_{(a,b)} \rangle < N) \leq P(\Omega \backslash E_{a,b}) \qquad \forall N > 0.$$

Letting $N \to \infty$, we get $P(\Omega \backslash E_{a,b}) = 1$ (because $\mathcal{L}_{\mathrm{stat}}$ is supported on $\mathcal{R}$) and then $P(E_{a,b}) = 0$. This implies that $P(\mathbf{U}_{\mathrm{stat}} \notin \mathcal{R}) = 0$.

Now, replacing $\mathcal{R}$ by $\mathcal{D}$ and $E_{a,b}$ by $E_{a,\infty}$, we obtain similarly that $P(\mathbf{U}_{\mathrm{stat}} \notin \mathcal{D}) = 0$ as soon as there exists a stationary distribution $\mathcal{L}_{\mathrm{stat}}$ such that $\mathcal{L}_{\mathrm{stat}}(\mathcal{D}) = 1$. $\square$

Let us now extend these results to Markov processes that take values in some $\sigma$-compact space $E$ and that do not necessarily satisfy the superposition property. In order to introduce some immigration and some superposition property, we will work on $\mathfrak{M}_E$, the set of point measures on $E$: if $\mathfrak{m} \in \mathfrak{M}_E$, either $\mathfrak{m} = \sum_{i \geq 1} \delta_{x^{(i)}}$ for some sequence $(x^{(i)}, i \geq 1)$ of points of $E$, or $\mathfrak{m} = \mathfrak{o}$, where $\mathfrak{o}$ is the trivial measure: $\mathfrak{o}(E) = 0$. The subset of measures of $\mathfrak{M}_E$ that are Radon is denoted by $\mathfrak{M}_E^{\text{Radon}}$ and is equipped with the topology of vague convergence. Consider then $I$, a nonnegative $\sigma$-finite measure on $E$, and $(X(t), t \geq 0)$, a Markov process with values in $E$. For any $\mathfrak{m} = \sum_{i \geq 1} \delta_{x^{(i)}} \in \mathfrak{M}_E$, set

$$\mathcal{X}^{(\mathfrak{m})}(t) := \sum_{i \geq 1} \delta_{X^{(x^{(i)})}(t)}, \qquad t \geq 0,$$

where $X^{(x^{(1)})}, X^{(x^{(2)})}, \ldots$ are independent versions of $X$, starting, respectively, from $x^{(1)}, x^{(2)}, \ldots$. If $\mathfrak{m} = \mathfrak{o}$, $\mathcal{X}^{(\mathfrak{m})}(t) := \mathfrak{o}$, $\forall t \geq 0$.

We now construct some $\mathcal{X}$-*process with immigration*. Let $\mathfrak{m}$ be a random element of $\mathfrak{M}_E^{\text{Radon}}$ and $((x(t_i), t_i), i \geq 1)$ be the atoms of a Poisson measure with intensity $I(d\mathbf{s})\, dt$, $t \geq 0$, independent of $\mathfrak{m}$. Conditionally on this Poisson measure and on $\mathfrak{m}$, let $\mathcal{X}^{(\mathfrak{m})}$ and $\mathcal{X}^{(\delta_{x(t_i)})}, i \geq 1$, be independent versions of $\mathcal{X}$ starting, respectively, from $\mathfrak{m}, \delta_{x(t_1)}, \delta_{x(t_2)}, \ldots$. Define then

$$\mathcal{X}\mathcal{I}^{(\mathfrak{m})}(t) := \mathcal{X}^{(\mathfrak{m})}(t) + \sum_{t_i \leq t} \mathcal{X}^{(\delta_{x(t_i)})}(t - t_i), \qquad t \geq 0,$$

and suppose that a.s., for all $t \geq 0$, $\mathcal{X}\mathcal{I}^{(\mathfrak{m})} \in \mathfrak{M}_E^{\text{Radon}}$. Then $\mathcal{X}\mathcal{I}^{(\mathfrak{m})}$ is Markovian and called $\mathcal{X}$-process with immigration starting from $\mathfrak{m}$. Introduce next the point measure

$$\mathcal{U}_{\text{stat}} := \sum_{t_i \geq 0} \mathcal{X}^{(\delta_{x(t_i)})}(t_i) = \sum_{i \geq 1} \delta_{X^{(x(t_i))}(t_i)}.$$

With the same kind of arguments as above, one obtains the following result.

LEMMA 6. (i) *Assume* $\mathcal{U}_{\text{stat}} \in \mathfrak{M}_E^{\text{Radon}}$ *a.s. Then the distribution* $\mathcal{L}(\mathcal{U}_{\text{stat}})$ *is a stationary distribution for* $\mathcal{X}\mathcal{I}$ *and* $\mathcal{X}\mathcal{I}^{(\mathfrak{m})}(t) \overset{\text{law}}{\to} \mathcal{U}_{\text{stat}}$ *as soon as* $\mathcal{X}^{(\mathfrak{m})}(t) \overset{\text{P}}{\to} \mathfrak{o}$ *as* $t \to \infty$.

(ii) *If* $P(\mathcal{U}_{\text{stat}} \notin \mathfrak{M}_E^{\text{Radon}}) > 0$, *there exists no stationary distribution for* $\mathcal{X}\mathcal{I}$.

2.2. *Conditions for existence and properties of FIs stationary distribution.* Up to now, $I$ is an immigration measure as defined in Section 1.2, that is, $I$ satisfies hypothesis (H1). Let *FI* denote a fragmentation with immigration $(\alpha, c, \nu, I)$. By definition, the fragmentation process satisfies the superposition property and, for every $\mathbf{u} \in \mathcal{D}$, $F^{(\mathbf{u})}(t) \overset{\text{a.s.}}{\to} \mathbf{0}$ as $t \to \infty$. Then the results of Lemma 5 can be rephrased as follows: if $((\mathbf{s}(t_i), t_i), i \geq 1)$ are the atoms of a Poisson measure with intensity $I(d\mathbf{s})\,dt$ and if conditionally on this Poisson measure, $F^{(\mathbf{s}(t_1))}, F^{(\mathbf{s}(t_2))}, \ldots$ are independent $(\alpha, c, \nu)$-fragmentations starting, respectively, from $\mathbf{s}(t_1), \mathbf{s}(t_2), \ldots$, then there is a stationary distribution for the fragmentation with immigration $(\alpha, c, \nu, I)$ if and only if

$$\mathbf{U}_{\text{stat}} = \sum_{t_i \geq 0} F^{(\mathbf{s}(t_i))}(t_i) \in \mathcal{D} \qquad \text{a.s.}$$

In this case,

$$FI^{(\mathbf{u})}(t) \stackrel{\text{law}}{\to} \mathbf{U}_{\text{stat}} \qquad \text{as } t \to \infty$$

for all $\mathbf{u} \in \mathcal{D}$ and, therefore, $\mathcal{L}(\mathbf{U}_{\text{stat}})$ is the *unique* stationary distribution for $FI$. The point is then to see when $\mathbf{U}_{\text{stat}}$ belongs to $\mathcal{D}$ and when it does not. The results are given in Section 2.2.1 where we further investigate whether $\mathbf{U}_{\text{stat}}$ is in $l^p$ or not, $p \geq 0$. This is particularly interesting when $\mathbf{U}_{\text{stat}} \in l^1$ a.s.: then the total mass of the system converges to an equilibrium, which means that the immigration compensates the mass lost by formation of dust (when $\alpha < 0$), by erosion or within sudden dislocations. When $\mathbf{U}_{\text{stat}} \in \mathcal{D}$ a.s., we also investigate the behavior of its small components. The proofs are detailed in Section 2.2.2.

2.2.1. *Statement of results.* Let $F$ denote a standard $(\alpha, c, \nu)$-fragmentation. In the statements below, we shall sometimes suppose that

(H2) $c = 0, \quad \nu\left(\sum_{j \geq 1} s_j < 1\right) = 0 \quad \text{and} \quad \int_{\mathcal{D}_1} \sum_{j \geq 1} |\ln(s_j)| s_j \nu(d\mathbf{s}) < \infty$

or

(H3) $\nexists\, 0 < r < 1 : F_i(t) \in \{r^n, n \in \mathbb{N}\} \quad \forall t \geq 0, i \geq 1$, and (H2) holds.

In terms of $\xi$, the subordinator driving a tagged fragment of $F$, the hypothesis (H2) means that $E[\xi(1)] < \infty$. We shall also use the convention $l^p = l^0$ when $p \leq 0$.

We now state our results on the existence of a stationary distribution; they depend heavily on the value of the index $\alpha$.

THEOREM 7. *Suppose $\alpha < 0$.*

(i) *If either $\int_{l^1} \sum_{j \geq 1} s_j^{-\alpha} \mathbf{1}_{\{s_j \geq 1\}} I(d\mathbf{s}) < \infty$ or $\int_{l^1} s_1^{-\alpha} \ln s_1 \mathbf{1}_{\{s_1 \geq 1\}} I(d\mathbf{s}) < \infty$, then the stationary state $\mathbf{U}_{\text{stat}} \in l^p$ a.s. for all $p > 1 + \alpha$.*
(ii) *There exists no stationary distribution when $\int_{l^1} s_1^{-\alpha} \mathbf{1}_{\{s_1 \geq 1\}} I(d\mathbf{s}) = \infty$.*

THEOREM 8. *Suppose $\alpha = 0$.*

(i) *If $\int_{l^1} \ln s_1 \mathbf{1}_{\{s_1 \geq 1\}} I(d\mathbf{s}) < \infty$, then, with probability one, $\mathbf{U}_{\text{stat}} \in l^p$ for all $p > 1$ and does not belong to $l^1$ when $c = 0$ and $\nu(\sum_{j \geq 1} s_j < 1) = 0$.*

(ii) *There exists no stationary distribution when $\int_{l^1} \ln s_1 \mathbf{1}_{\{s_1 \geq 1\}} I(d\mathbf{s}) = \infty$ and* (H2) *holds.*

THEOREM 9.  *Suppose $\alpha > 0$. If $\int_{l^1} s_1^\varepsilon \mathbf{1}_{\{s_1 \geq 1\}} I(d\mathbf{s}) < \infty$ for some $\varepsilon > 0$, then $\mathbf{U}_{\text{stat}} \in l^p$ a.s. for $p$ large enough and if* (H3) *holds, then $\mathbf{U}_{\text{stat}} \notin l^{1+\alpha}$ a.s. More precisely, for every $\gamma > 0$:*

(i) *if $\int_{l^1} \sum_{j \geq 1} s_j^\gamma \mathbf{1}_{\{s_j \geq 1\}} I(d\mathbf{s}) < \infty$, then $\mathbf{U}_{\text{stat}} \in l^p$ a.s. for all $p > 1 + \alpha/(\gamma \wedge 1)$,*

(ii) *if $\int_{l^1} s_1^\gamma \mathbf{1}_{\{s_1 \geq 1\}} I(d\mathbf{s}) = \infty$ and* (H3) *holds, then $\mathbf{U}_{\text{stat}} \notin l^{1+\alpha/(\gamma \wedge 1)}$ a.s.*

When $-1 < \alpha < 0$, the result of Theorem 7(i) can be completed (see the remark following Proposition 10 below): in most cases, either $\mathbf{U}_{\text{stat}} \notin l^{1+\alpha}$ a.s. or both events $\{\mathbf{U}_{\text{stat}} = \mathbf{0}\}$ and $\{\mathbf{U}_{\text{stat}} \notin l^{1+\alpha}\}$ have positive probabilities.

It is interesting to notice that the above conditions for existence or absence of a stationary distribution depend only on $\alpha$ and $I$, provided hypothesis (H3) holds. Indeed, recall the definition (1) of $\alpha_I$ and let then $\alpha$ vary. According to the above theorems, the values $\alpha = \alpha_I$ and $\alpha = -1$ are critical. Provided $\alpha_I < 0$, the stationary distribution exists when $\alpha > \alpha_I$ and does not exist when $\alpha < \alpha_I$. Moreover, the stationary state $\mathbf{U}_{\text{stat}}$ is a.s. composed by a finite number of particles as soon as $\alpha_I < \alpha < -1$, whereas when $\alpha > -1$, $\mathbf{U}_{\text{stat}} \notin l^{1+\alpha}$ with a positive probability (which equals 1 when $\alpha \geq 0$ and depends on further hypothesis on $I$ and $\alpha$ when $-1 < \alpha < 0$, see the forthcoming Proposition 10).

Let us try to explain these results. By the scaling property of fragmentation processes, particles with mass $\geq 1$ split faster when $\alpha$ is larger. This explains that, when $\alpha$ is too small, some particles may accumulate in intervals of type $(a, \infty)$, $a > 0$, which implies that $\mathbf{U}_{\text{stat}} \notin \mathcal{D}$. For $\alpha$ large enough, particles with mass $\geq 1$ become rapidly smaller, but particles with mass $\leq 1$ split more slowly when $\alpha$ is larger. Therefore, small particles accumulate and $\mathbf{U}_{\text{stat}} \notin l^p$ when $p$ is too small. Moreover, the smallest $p$ such that $\mathbf{U}_{\text{stat}} \in l^p$ increases as $\alpha$ increases. When $\alpha < -1$, it is known that small particles are very quickly reduced to dust (see, e.g., Proposition 2 in [8]). This implies that $\mathbf{U}_{\text{stat}} \in l^0$, provided it belongs to $\mathcal{D}$.

*Small particles behavior.* Suppose that $-1 < \alpha < 0$ and $\int_{l^1} \sum_{j \geq 1} s_j^{-\alpha} \mathbf{1}_{\{s_j \geq 1\}} \times I(d\mathbf{s}) < \infty$, so that $\mathbf{U}_{\text{stat}} \in \mathcal{D}$ a.s., according to Theorem 7(i). Consider then the random function

$$\varepsilon \mapsto \overline{\mathbf{U}}_{\text{stat}}(\varepsilon) := \mathbf{U}_{\text{stat}}([\varepsilon, \infty)),$$

which counts the number of components of $\mathbf{U}_{\text{stat}}$ larger than $\varepsilon$. We want to investigate the limiting behavior of $\overline{\mathbf{U}}_{\text{stat}}(\varepsilon)$ as $\varepsilon \to 0$. In that aim, we make the following technical hypothesis:

(H4)
$$\int_{\mathcal{D}_1} \sum_{j>i\geq 1} s_i^{1+\alpha} s_j \nu(d\mathbf{s}) < \infty \quad \text{and}$$
$$\int_{\mathcal{D}_1} (1-s_1)^\theta \nu(d\mathbf{s}) < \infty \qquad \text{for some } \theta < 1,$$

as well as hypothesis (H3). Note that the first integral involved in (H4) is finite as soon as $\alpha > -1$ and $\nu(s_N > 0) = 0$ for some integer $N \geq 2$, because then $\int_{\mathcal{D}_1} \sum_{j>i\geq 1} s_i^{1+\alpha} s_j \nu(d\mathbf{s}) \leq (N-1) \int_{\mathcal{D}_1} (1-s_1) \nu(d\mathbf{s})$.

PROPOSITION 10. *Under the previous hypotheses:*

(i) *if $\int_{l^1} \sum_{j\geq 1} s_j^{-\alpha} \mathbf{1}_{\{s_j \leq 1\}} I(d\mathbf{s}) < \infty$, there exists a finite r.v. $X$, $0 < P(X=0) < 1$, such that*
$$\overline{\mathbf{U}}_{\text{stat}}(\varepsilon) \varepsilon^{1+\alpha} \underset{\varepsilon \to 0}{\to} X \qquad a.s.,$$

(ii) *if $\int_{l^1} s_1^{-\alpha} \mathbf{1}_{\{s_j \leq 1\}} I(d\mathbf{s}) = \infty$, one has $\liminf_{\varepsilon \to 0} \varepsilon^{1+\alpha} \overline{\mathbf{U}}_{\text{stat}}(\varepsilon) > 0$ a.s.*

In particular, this implies that $P(\mathbf{U}_{\text{stat}} \notin l^{1+\alpha}) = 1$ when the assumption of the second statement is satisfied. This is not true when the assumption of the first statement holds: in such case, $0 < P(\mathbf{U}_{\text{stat}} = \mathbf{0}) \leq P(\mathbf{U}_{\text{stat}} \in l^{1+\alpha}) < 1$ [see the proof of (i) for the first inequality].

When $\alpha \geq 0$ or $\alpha < -1$, some information on the behavior of $\overline{\mathbf{U}}_{\text{stat}}(\varepsilon)$ as $\varepsilon \to 0$ can be deduced from Theorems 7, 8 and 9. Thus, $\overline{\mathbf{U}}_{\text{stat}}(0) < \infty$ a.s. as soon as $\alpha_I < \alpha < -1$. To obtain some information when $\alpha \geq 0$, first notice that, when $\mathbf{U}_{\text{stat}} \in \mathcal{D}$, $\int_{(0,\infty)} x^p \mathbf{U}_{\text{stat}}(dx) < \infty \Leftrightarrow \int_{(0,1)} \overline{\mathbf{U}}_{\text{stat}}(x^{1/p}) dx < \infty$, by integration by parts. Combined with Theorem 9, this implies, when $\alpha > 0$, that if $\int_{l^1} \sum_{j\geq 1} s_j^\gamma \mathbf{1}_{\{s_j \geq 1\}} I(d\mathbf{s}) < \infty$, then $\liminf_{\varepsilon \to 0} \varepsilon^p \overline{\mathbf{U}}_{\text{stat}}(\varepsilon) = 0$ for all $p > 1 + \alpha/(\gamma \wedge 1)$, whereas if $\int_{l^1} s_1^\gamma \mathbf{1}_{\{s_1 \geq 1\}} I(d\mathbf{s}) = \infty$, $\limsup_{\varepsilon \to 0} \varepsilon^p \overline{\mathbf{U}}_{\text{stat}}(\varepsilon) = \infty$ for all $p < 1 + \alpha/(\gamma \wedge 1)$. The behavior near 0 of $\overline{\mathbf{U}}_{\text{stat}}(\varepsilon)$ is then strongly connected to the immigration $I$. Similarly, when $\alpha = 0$ and when there is a stationary distribution, one deduces from Theorem 8 that $\liminf_{\varepsilon \to 0} \varepsilon^p \overline{\mathbf{U}}_{\text{stat}}(\varepsilon) = 0$ for all $p > 1$, and that $\limsup_{\varepsilon \to 0} \varepsilon^p \overline{\mathbf{U}}_{\text{stat}}(\varepsilon) = \infty$ for all $p < 1$, provided $c = \nu(\sum_{i \geq 1} s_i < 1) = 0$.

REMARK. It is possible to show that $\mathbf{U}_{\text{stat}} \in \mathcal{R}$ a.s. as soon as $\int_{l^1} \sum_{j\geq 1} s_j \times \mathbf{1}_{\{s_j \geq 1\}} I(d\mathbf{s}) < \infty$ and that $P(\mathbf{U}_{\text{stat}} \notin \mathcal{R}) > 0$ as soon as $\alpha > -1$, $\int_{l^1} s_1^{-\alpha} \mathbf{1}_{\{s_1 \geq 1\}} \times I(d\mathbf{s}) = \infty$ and hypotheses (H3) and (H4) hold. The first claim can be proved by using some arguments of the proof of the forthcoming Proposition 16 and the second claim is a consequence of Theorems 4(i) and 7 of [19], which are also used below to prove Proposition 10.

2.2.2. *Proofs.* Let $F$ be a standard $(\alpha, c, \nu)$-fragmentation and for every $p \in \mathbb{R}$ and $t \geq 0$, define

$$M(p,t) := \sum_{k \geq 1} (F_k(t))^p \mathbf{1}_{\{F_k(t) > 0\}},$$

which is a.s. finite at least when $p \geq 1$ (since it is bounded from above by 1). That $\mathbf{U}_{\text{stat}}$ belongs to some $l^p$-space is closely related to the behavior of the function $t \mapsto M(p,t)$. Indeed,

$$\mathbf{U}_{\text{stat}} = \sum_{i \geq 1} \sum_{j \geq 1} s_j(t_i) F^{(i,j)}(s_j^\alpha(t_i)t_i),$$

where the $F^{(i,j)}$'s, $i,j \geq 1$, are i.i.d. copies of $F$, independent of $((\mathbf{s}(t_i), t_i), i \geq 1)$. Then $\mathbf{U}_{\text{stat}} \in l^p \Leftrightarrow M(p) < \infty$ with

$$\begin{aligned} M(p) &= \int_{(0,\infty)} x^p \mathbf{U}_{\text{stat}}(dx) \\ &= \sum_{i \geq 1} \sum_{j \geq 1} s_j^p(t_i) M^{(i,j)}(p, s_j^\alpha(t_i)t_i) \mathbf{1}_{\{s_j(t_i) > 0\}}, \end{aligned}$$

where the $M^{(i,j)}(p, \cdot)$'s, $i,j \geq 1$, are i.i.d. copies of $M(p, \cdot)$, independent of $((\mathbf{s}(t_i), t_i), i \geq 1)$. Using the tagged particle approach as explained in Section 1.1, one obtains the following results on $M(p, \cdot)$.

LEMMA 11. (i) *Suppose* $\alpha \leq 0$. *Then* $\int_0^\infty \exp(\lambda t) E[M(p,t)] \, dt < \infty$ *as soon as* $p \geq 1 + \alpha$ *and* $\lambda < \phi(p - 1 - \alpha)$. *In particular,* $E[M(p,t)] < \infty$ *for a.e.* $t \geq 0$ *as soon as* $p \geq 1 + \alpha$.

(ii) *Suppose* $\alpha > 0$. *Then for every* $\eta > 0$ *and every* $p \geq 1$, *there exists a random variable* $D_{(\eta,p)}$ *with positive moments of all orders such that*

$$M(p,t) \leq D_{(\eta,p)} t^{-(p-1)/(\alpha+\eta)} \qquad \text{a.s. for every } t > 0.$$

*Consequently,* $\int_0^\infty E[M(p,t)] \, dt < \infty$ *when* $p > 1 + \alpha$.

Bertoin (Corollary 3 in [8]) shows that when $\alpha > 0$ and $p \geq 1$, the process $t^{(p-1)/\alpha} M(p,t)$ converges in probability to some deterministic limit as $t \to \infty$, provided the fragmentation satisfies hypothesis (H3). See also Brennan and Durrett [11, 12], who prove the almost sure convergence for binary fragmentations ($\nu(s_1 + s_2 < 1) = 0$) with a finite dislocation measure.

PROOF. We use the notation introduced in Section 1.1.
(i) According to (6),

$$E[M(p,t)] = E[\exp((1-p)\overline{\xi}(\rho(t)))\mathbf{1}_{\{t<D\}}],$$

where $D = \inf\{t : \rho(t) \geq \mathbf{e}(k)\}$. Therefore,

$$\int_0^\infty \exp(\lambda t) E[M(p,t)] \, dt$$

(14)
$$= E\left[\int_0^D \exp(\lambda t) \exp((1-p)\overline{\xi}(\rho(t))) \, dt\right]$$

$$= E\left[\int_0^{\mathbf{e}(k)} \exp(\lambda \rho^{-1}(t)) \exp((1-p+\alpha)\overline{\xi}(t)) \, dt\right],$$

using for the last equality the change of variables $t \mapsto \rho(t)$ and that, by definition of $\rho$, $\exp(\alpha\overline{\xi}(\rho(t))) \, d\rho(t) = dt$ on $[0, D)$. The function $\rho^{-1}$ denotes the right inverse of $\rho$ and, clearly, $\rho^{-1}(t) \leq t$ since $\alpha \leq 0$. When $p \geq 1 + \alpha$, this leads to

$$\int_0^\infty \exp(\lambda t) E[M(p,t)] \, dt$$

$$\leq \begin{cases} E\left[\int_0^{\mathbf{e}(k)} \exp(-\overline{\phi}(p-1-\alpha)t) \, dt\right], & \text{if } \lambda < 0, \\ E\left[\int_0^{\mathbf{e}(k)} \exp((\lambda - \overline{\phi}(p-1-\alpha))t) \, dt\right], & \text{if } \lambda \geq 0, \end{cases}$$

and in both cases, the integral is finite as soon as $\lambda < \phi(p-1-\alpha) = \overline{\phi}(p-1-\alpha) + k$.

(ii) Fix $\alpha > 0$, $p \geq 1$ and $\eta > 0$ and recall that, according to (5),

$$M(p,t) = E[\exp(-(p-1)\xi(\rho(t)))\mathbf{1}_{\{t<D\}}|F].$$

Since $\xi$ is increasing, one has

$$\rho(t) \exp(-\eta\xi(\rho(t))) \leq \int_0^{\rho(t)} \exp(-\eta\xi(r)) \, dr \leq \int_0^\infty \exp(-\eta\xi(r)) \, dr := D_{(\eta)}.$$

And, on the other hand, for $t < D$,

$$t = \int_0^{\rho(t)} \exp(\alpha\xi(r)) \, dr \leq \rho(t)\exp(\alpha\xi(\rho(t))).$$

Combining these inequalities, we obtain $\exp(-(\alpha+\eta)\xi(\rho(t))) \leq t^{-1}D_{(\eta)}$ for all $t < D$. Hence, $M(p,t) \leq t^{-(p-1)/(\alpha+\eta)}D_{(\eta,p)}$, where $D_{(\eta,p)} := E[D_{(\eta)}^{(p-1)/(\alpha+\eta)}|F]$. Carmona, Petit and Yor [13] have shown that $D_{(\eta)}$ has moments of all positive orders, which, by Hölder's inequality, is also true for $D_{(\eta,p)}$. □

We now turn to the proofs of Theorems 7, 8 and 9.

PROOF OF THEOREM 7. (i) Fix $p > 1 + \alpha$ and split $M(p)$ into two subsums:

$$M_{\inf}(p) = \sum_{i \geq 1}\sum_{j \geq 1} s_j^p(t_i)\mathbf{1}_{\{0<s_j(t_i)<1\}} M^{(i,j)}(p, s_j^\alpha(t_i)t_i)$$

and $M_{\sup}(p) = M(p) - M_{\inf}(p)$. One has

$$E[M_{\inf}(p)] = \int_{l^1} \left(\sum_{j \geq 1} s_j^{p-\alpha} \mathbf{1}_{\{s_j < 1\}}\right) I(d\mathbf{s}) \times \int_0^\infty E[M(p,t)]\,dt$$

and both of these integrals are finite according to hypothesis (H1) and Lemma 11 since $p > 1 + \alpha$. It remains to show that $M_{\sup}(p) < \infty$ when $I$ integrates $\sum_{j \geq 1} s_j^{-\alpha} \mathbf{1}_{\{s_j \geq 1\}}$ or $s_1^{-\alpha} \ln s_1 \mathbf{1}_{\{s_1 \geq 1\}}$.

Suppose first that $\int_{l^1} \sum_{j \geq 1} s_j^{-\alpha} \mathbf{1}_{\{s_j \geq 1\}} I(d\mathbf{s}) < \infty$ and let $\tau^{(i,j)}$ be the first time at which the fragmentation $F^{(i,j)}$ is entirely reduced to dust. Equivalently, $\tau^{(i,j)}$ is the first time at which $M^{(i,j)}$ reaches 0. If the number of pairs $(i,j)$ such that $s_j^\alpha(t_i)t_i \leq \tau^{(i,j)}$ and $s_j(t_i) \geq 1$ is finite, then the sum $M_{\sup}(p)$ is finite because it involves at most a finite number of nonzero $M^{(i,j)}(p, s_j^\alpha(t_i)t_i)$ [which are a.s. all finite according to Lemma 11(i)]. To prove that this is the case, we use the theory of Poisson measures. Since the r.v. $\tau^{(i,j)}$, $i,j \geq 1$, are i.i.d., the measure

$$\sum_{i \geq 1} \delta_{t_i^{-1} \sup_{j\,:\,s_j(t_i) \geq 1}(\tau^{(i,j)} s_j^{-\alpha}(t_i))}$$

is a Poisson measure with intensity $m$ defined for any positive measurable function $f$ by

$$\int_0^\infty f(x) m(dx) = \int_0^\infty \int_{l^1} E\left[f\left(t^{-1} \sup_{j\,:\,s_j \geq 1}(\tau^{(1,j)} s_j^{-\alpha})\right)\right] I(d\mathbf{s})\,dt.$$

The integral $\int_1^\infty m(dx)$ is bounded from above by $E[\tau^{(1,1)}] \int_{l^1} \sum_{j \geq 1} s_j^{-\alpha} \mathbf{1}_{\{s_j \geq 1\}} \times I(d\mathbf{s})$, which is finite by assumption on $I$ and since $E[\tau^{(1,1)}] < \infty$ [by (8)]. This implies that a.s. there are only a finite number of integers $i \geq 1$ such that $t_i^{-1} \sup_{j\,:\,s_j(t_i) \geq 1}(\tau^{(i,j)} s_j^{-\alpha}(t_i)) \geq 1$. For each of these $i$, there is at most a finite number of integers $j \geq 1$ such that $s_j(t_i) \geq 1$. Hence, the number of pairs $(i,j)$ such that $s_j^\alpha(t_i)t_i \leq \tau^{(i,j)}$ and $s_j(t_i) \geq 1$ is indeed a.s. finite.

Assume now that $\int_{l^1} s_1^{-\alpha} \ln s_1 \mathbf{1}_{\{s_1 \geq 1\}} I(d\mathbf{s}) < \infty$. For any $a > 0$, the number of integers $i \geq 1$ such that $at_i \leq s_1^{-\alpha}(t_i) \ln(s_1(t_i))$ and $s_1(t_i) \geq 1$ is then a.s. finite. The sum $M_{\sup}(p)$ is therefore finite if

$$\sum_{i \geq 1} \sum_{j \geq 1} s_j^p(t_i) \mathbf{1}_{\{at_i > s_1^{-\alpha}(t_i) \ln(s_1(t_i))\}} \mathbf{1}_{\{s_j(t_i) \geq 1\}} M^{(i,j)}(p, s_j^\alpha(t_i) t_i)$$

is finite for some (and then all) $a > 0$. The expectation of this latter sum is bounded from above by

$$\int_0^\infty \int_{l^1} \left(\sum_{j \geq 1} s_j^p \mathbf{1}_{\{at > s_j^{-\alpha} \ln s_j\}} \mathbf{1}_{\{s_j \geq 1\}}\right) E[M(p, s_j^\alpha t)] I(d\mathbf{s})\,dt \qquad (\text{as } s_j \leq s_1)$$

$$\leq \int_{l^1} \sum_{j \geq 1} \mathbf{1}_{\{s_j \geq 1\}} I(d\mathbf{s}) \int_0^\infty \exp(at(p-\alpha)) E[M(p,t)]\,dt,$$

which is finite for $a$ sufficiently small, according to Lemma 11(i). Hence, $M_{\sup}(p) < \infty$ a.s.

(ii) Suppose $\int_{l^1} s_1^{-\alpha} \mathbf{1}_{\{s_1 \geq 1\}} I(d\mathbf{s}) = \infty$ and let $\tau_{1/2}^{(i,1)} := \inf\{t \geq 0 : F_1^{(i,1)}(t) < 1/2\}$ be the first time at which all components of $F^{(i,1)}$ are smaller than $1/2$, $i \geq 1$. Note that $E[\tau_{1/2}^{(i,1)}] > 0$ since $F_1^{(i,1)}$ is càdlàg. The measure

$$\sum_{i \geq 1 \,:\, s_1(t_i) \geq 1} \delta_{s_1^{-\alpha}(t_i) t_i^{-1} \tau_{1/2}^{(i,1)}}$$

is a Poisson measure with intensity $m'$ given by

$$\int_0^\infty f(x) m'(dx) = \int_0^\infty \int_{l^1} E[f(s_1^{-\alpha} t^{-1} \tau_{1/2}^{(1,1)})] \mathbf{1}_{\{s_1 \geq 1\}} I(d\mathbf{s}) \, dt.$$

By assumption on $I$ and since $E[\tau_{1/2}^{(1,1)}] > 0$, the integral $\int_1^\infty m'(dx)$ is infinite and, consequently, the number of integers $i$ such that $\tau_{1/2}^{(i,1)} > s_1^\alpha(t_i) t_i$ and $s_1(t_i) \geq 1$ is a.s. infinite. For those $i$, $s_1(t_i) F_1^{(i,1)}(s_1^\alpha(t_i) t_i) \geq 1/2$ and, therefore, $\mathbf{U}_{\mathrm{stat}}$ contains a sequence of terms all larger than $1/2$, which implies that it is not in $\mathcal{D}$ a.s. $\square$

PROOF OF THEOREM 8. (i) The second part of the proof of Theorem 7(i) (replacing there $\alpha$ by 0) shows that $\mathbf{U}_{\mathrm{stat}} \in \bigcap_{p>1} l^p$ when $\int_{l^1} \ln(s_1) \mathbf{1}_{\{s_1 \geq 1\}} I(d\mathbf{s}) < \infty$. Now, if $c = 0$ and $\nu(\sum_{k \geq 1} s_k < 1) = 0$, the sum $M(1)$ equals $\sum_{i \geq 1} \sum_{j \geq 1} s_j(t_i)$, which is clearly a.s. infinite since $I \neq 0$.

(ii) Assume that $\int_{l^1} \ln(s_1) \mathbf{1}_{\{s_1 \geq 1\}} I(d\mathbf{s}) = \infty$ and $E[\xi(1)] < \infty$. For each $i \geq 1$, let $\exp(-\xi^{(i,1)}(\cdot))$ denote the process of masses of the tagged particle in the fragmentation $F^{(i,1)}$. To prove that $\mathbf{U}_{\mathrm{stat}} \notin \mathcal{D}$, it suffices to show that its subsequence $\{s_1(t_i) \exp(-\xi^{(i,1)}(t_i)), i \geq 1\}^\downarrow \notin \mathcal{D}$. The components of this sequence are the atoms of a Poisson measure with intensity $m''$ given by

$$\int_0^\infty f(x) m''(dx) = \int_0^\infty \int_{l^1} E[f(s_1 \exp(-\xi(t)))] I(d\mathbf{s}) \, dt.$$

Take then $a > E[\xi(1)]$. Since $\xi(t)/t \stackrel{\mathrm{a.s.}}{\to} E[\xi(1)]$ as $t \to \infty$, there exists some $t_0$ such that $P(\xi(t) \leq at) \geq 1/2$ for $t \geq t_0$. Then

$$\int_1^\infty m''(dx) = \int_0^\infty \int_{l^1} P(\xi(t) \leq \ln s_1) I(d\mathbf{s}) \, dt$$

$$\geq \int_{l^1} \int_{t_0}^{a^{-1} \ln s_1} P(\xi(t) \leq at) \, dt \, I(d\mathbf{s})$$

$$\geq \tfrac{1}{2} \int_{l^1} (a^{-1} \ln s_1 - t_0) \mathbf{1}_{\{a^{-1} \ln s_1 \geq t_0\}} I(d\mathbf{s})$$

and this last integral is infinite by assumption. Hence, $\sum_{i\geq 1} \delta_{s_1(t_i)\exp(-\xi^{(i,1)}(t_i))} \notin \mathcal{D}$ a.s. and a fortiori $\mathbf{U}_{\text{stat}} \notin \mathcal{D}$ a.s. □

PROOF OF THEOREM 9. Fix $p \geq 1 + \alpha$. According to the Campbell formula for Poisson measures (see [21]), the sum $M(p)$ is finite if and only if

$$\int_0^\infty \int_{l^1} E\left[1 - \exp\left(-\sum_{j\geq 1} s_j^p M^{(1,j)}(p, s_j^\alpha t)\right)\right] I(d\mathbf{s})\, dt < \infty. \tag{15}$$

(i) We first prove assertion (i) and that $\mathbf{U}_{\text{stat}} \in l^p$ a.s. for $p$ large enough when $I$ integrates $s_1^\varepsilon \mathbf{1}_{\{s_1 \geq 1\}}$. Suppose $p > 1 + \alpha$ and note that the integral (15) is bounded from above by

$$\int_{l^1} \sum_{j\geq 1} s_j^{p-\alpha} \mathbf{1}_{\{s_j < 1\}} I(d\mathbf{s}) \int_0^\infty E[M(p, t)]\, dt$$

$$+ \int_0^\infty \int_{l^1} E\left[1 - \exp\left(-\sum_{j\geq 1} s_j^p \mathbf{1}_{\{s_j \geq 1\}} M^{(1,j)}(p, s_j^\alpha t)\right)\right] I(d\mathbf{s})\, dt.$$

According to Lemma 11(ii), the first component of this sum is finite and, for all $\eta > 0$, there exists some i.i.d. r.v. $D_{(\eta,p)}^{(j)}$ having finite moments of all positive orders and independent of $(\mathbf{s}(t_i), i \geq 1)$ such that the second component is bounded from above by

$$\int_0^\infty \int_{l^1} E\left[1 - \exp\left(-\sum_{j\geq 1} s_j^{p-\alpha(p-1)/(\alpha+\eta)} \mathbf{1}_{\{s_j\geq 1\}} D_{(\eta,p)}^{(j)} t^{-(p-1)/(\alpha+\eta)}\right)\right] I(d\mathbf{s})\, dt.$$

Using the fact that $\int_0^\infty E[1-\exp(-t^{-a}X)]\, dt = \int_0^\infty (1-\exp(-t^{-a}))\, dt\, E[X^{1/a}]$ for nonnegative r.v. $X$, one sees that this double integral is equal to

$$\int_0^\infty (1 - \exp(-t^{-(p-1)/(\alpha+\eta)}))\, dt$$

$$\times \int_{l^1} \left(\sum_{j\geq 1} s_j^{(p\eta+\alpha)/(\alpha+\eta)} \mathbf{1}_{\{s_j\geq 1\}}\right)^{(\alpha+\eta)/(p-1)} I(d\mathbf{s}) E[D_{(\eta,p)}^{(1)}{}^{(\alpha+\eta)/(p-1)}].$$

If $p > 1 + \alpha + \eta$, the first integral in this latter product is finite. So, take $\eta$ small enough so that $p > 1 + \alpha + \eta$ and notice then that

$$\int_{l^1} \left(\sum_{j\geq 1} s_j^{(p\eta+\alpha)/(\alpha+\eta)} \mathbf{1}_{\{s_j\geq 1\}}\right)^{(\alpha+\eta)/(p-1)} I(d\mathbf{s}) \tag{16}$$

$$\leq \int_{l^1} \sum_{j\geq 1} s_j^{(p\eta+\alpha)/(p-1)} \mathbf{1}_{\{s_j\geq 1\}} I(d\mathbf{s}).$$

The integral (15) is therefore finite as soon as the integral in the right-hand side of (16) is finite for some $\eta > 0$ small enough. Hence, we get (i).

The same argument shows that $\mathbf{U}_{\text{stat}} \in l^p$ for $p$ sufficiently large when there exists some $\varepsilon > 0$ such that $\int_{l^1} s_1^\varepsilon \mathbf{1}_{\{s_1 \geq 1\}} I(d\mathbf{s}) < \infty$. Indeed, let $p > 1 + \alpha + \eta$. It suffices then to show that the integral on the left-hand of (16) is finite and to do so, we replace the upper bound there by

$$\int_{l^1} \left( \sum_{j \geq 1} s_j^{(p\eta+\alpha)/(\alpha+\eta)} \mathbf{1}_{\{s_j \geq 1\}} \right)^{(\alpha+\eta)/(p-1)} I(d\mathbf{s})$$

$$\leq \int_{l^1} s_1^{(p\eta+\alpha)(p-1)} \left( \sum_{j \geq 1} \mathbf{1}_{\{s_j \geq 1\}} \right)^{(\alpha+\eta)/(p-1)} I(d\mathbf{s}),$$

which, by Hölder's inequality, is finite as soon as $p$ is large enough and $\eta$ small enough.

(ii) We now turn to the proof of assertion (ii) and that $\mathbf{U}_{\text{stat}} \notin l^{1+\alpha}$ when (H3) holds. The integral (15) is bounded from below by

$$\int_0^\infty \int_{l^1} s_1^{-\alpha} E[(1 - \exp(-s_1^p M(p,t))) \mathbf{1}_{\{M(p,t) \geq rt^{-(p-1)/\alpha}\}}] I(d\mathbf{s}) \, dt$$

$$\geq \int_{l^1} s_1^{-\alpha} \int_0^\infty (1 - \exp(-s_1^p r t^{-(p-1)/\alpha})) P(M(p,t) \geq rt^{-(p-1)/\alpha}) \, dt \, I(d\mathbf{s}).$$

According to Corollary 3 in [8], the hypothesis (H3) ensures that $t^{(p-1)/\alpha} M(p,t)$ converges in probability to some finite deterministic constant as $t \to \infty$. Hence, taking $r > 0$ small enough and then $t_0$ large enough, one has $P(M(p,t) \geq rt^{-(p-1)/\alpha}) \geq 1/2$ for $t \geq t_0$ and, therefore, the integral (15) is bounded from below by

$$\tfrac{1}{2} \int_{l^1} s_1^{-\alpha} s_1^{p\alpha/(p-1)} \int_0^\infty \mathbf{1}_{\{s_1^{p\alpha/(p-1)} \geq (t_0/t)\}} (1 - \exp(-rt^{-(p-1)/\alpha})) \, dt \, I(d\mathbf{s}),$$

which is infinite as soon as $p \leq 1 + \alpha$ or $\int_{l^1} s_1^{\alpha/(p-1)} \mathbf{1}_{\{s_1 \geq t_0\}} I(d\mathbf{s}) = \infty$. □

PROOF OF PROPOSITION 10. For the standard fragmentation $F$, let $N_{(\varepsilon,\infty)}(t) := \sum_{k \geq 1} \mathbf{1}_{\{F_k(t) > \varepsilon\}}$ denote the number of terms larger than $\varepsilon$ present at time $t$. Under the hypotheses (H3), (H4) and $\alpha > -1$, Theorems 4(i) and 7 of [19] describe the behavior of $N_{(\varepsilon,\infty)}(t)$ as $\varepsilon \to 0$. Theorem 4(i) states the existence of a random function $L$ such that $\sum_{k \geq 1} F_k(t) = \int_t^\infty L(u) \, du$ a.s. for all $t$. Then Theorem 7 says that

(17) $$\varepsilon^{1+\alpha} N_{(\varepsilon,\infty)}(t) \to KL(t) \qquad \text{as } \varepsilon \to 0$$

a.s. for almost every $t$, where $K = (1+\alpha)/\alpha^2 E[\xi(1)]$. Note that the sum $\overline{\mathbf{U}}_{\text{stat}}(\varepsilon)$ can be written as

(18) $$\overline{\mathbf{U}}_{\text{stat}}(\varepsilon) = \sum_{i,j \geq 1} N^{(i,j)}_{(\varepsilon/s_j(t_i),\infty)}(s_j^\alpha(t_i) t_i),$$

where the $N^{(i,j)}_{(\cdot,\infty)}(\cdot)$'s are i.i.d. copies of $N_{(\cdot,\infty)}(\cdot)$, independent of $((\mathbf{s}(t_i),t_i), i \geq 1)$.

(i) Let $\tau^{(i,j)}$ be the first time at which $F^{(i,j)}$ reaches $\mathbf{0}$, $i,j \geq 1$. With the same arguments as in the proof of Theorem 7(i), one sees that, with probability one, there is at most a finite number of $t_i < \sup_{j \geq 1}(\tau^{(i,j)} s_j^{-\alpha}(t_i))$ if and only if $\int_{l^1} E[\sup_{j \geq 1} \tau^{(1,j)} s_j^{-\alpha}] I(d\mathbf{s}) < \infty$. This integral is finite by assumption. A moment of thought then shows that there is at most a finite number of integers $i,j \geq 1$—independent of $\varepsilon$—such that $N^{(i,j)}_{(\varepsilon/s_j(t_i),\infty)}(s_j^\alpha(t_i)t_i) > 0$. Consequently, the sum (18) involves a finite number of nonzero terms and

$$\varepsilon^{1+\alpha}\overline{\mathbf{U}}_{\mathrm{stat}}(\varepsilon) \underset{\varepsilon \to 0}{\to} K \sum_{i,j \geq 1} L^{(i,j)}(s_j^\alpha(t_i)t_i) s_j^{1+\alpha}(t_i) \qquad \text{a.s.,}$$

where the functions $L^{(i,j)}$'s are i.i.d. and distributed as $L$. This limit, which we denote by $X$, is null as soon as $\mathbf{U}_{\mathrm{stat}} = \mathbf{0}$, that is, as soon as there is no integer $i \geq 1$ such that $t_i < \sup_{j \geq 1}(\tau^{(i,j)} s_j^{-\alpha}(t_i))$. This occurs, according to the Poissonian construction, with a positive probability. On the other hand, the Lebesgue measure of $\mathcal{B}_L := \{x \geq 0 : L(x) > 0\}$ [denoted by $\mathrm{Leb}(\mathcal{B}_L)$] is a.s. nonzero and then $P(X > 0) > 0$.

(ii) Suppose $\int_{l^1} s_1^{-\alpha} \mathbf{1}_{\{s_1 \leq 1\}} I(d\mathbf{s}) = \infty$ and let $\mathcal{B}_{L^{(i,j)}} := \{x \geq 0 : L^{(i,j)}(x) > 0\}$, which are i.i.d. copies of $\mathcal{B}_L$. One checks that there a.s. exists a time $t_i \in \bigcup_{j \geq 1} s_j^{-\alpha}(t_i) \mathcal{B}_{L^{(i,j)}}$ if and only if the integral $\int_{l^1} E[\mathrm{Leb}(\bigcup_{j \geq 1} s_j^{-\alpha} \mathcal{B}_{L^{(1,j)}})] I(d\mathbf{s})$ is infinite and that this integral is indeed infinite here, according to the assumption on $I$ and since $\mathrm{Leb}(\mathcal{B}_L) > 0$ a.s. From this, we deduce that

$$\sum_{1 \leq i,j \leq N} L^{(i,j)}(s_j^\alpha(t_i)t_i) s_j^{1+\alpha}(t_i) > 0 \qquad \text{a.s. for } N \text{ large enough}$$

and then, by (17) and (18), that $\liminf_{\varepsilon \to 0} \varepsilon^{1+\alpha} \overline{\mathbf{U}}_{\mathrm{stat}}(\varepsilon) > 0$. $\square$

**3. Rate of convergence to the stationary distribution.** We are interested in the convergence in law to the stationary regime $\mathbf{U}_{\mathrm{stat}}$. It is already known, according to Lemma 5, that, for every random $\mathbf{u} \in \mathcal{D}$, the process $FI^{(\mathbf{u})}(t)$ converges in law as $t \to \infty$ to the stationary state $\mathbf{U}_{\mathrm{stat}}$, provided it belongs to $\mathcal{D}$ a.s. The aim of this section is to strengthen this result by providing upper bounds for the rate at which this convergence takes place. The norm considered on the set of signed finite measures on $\mathcal{D}$ is

$$\|\mu\| := \sup_{f \text{1-Lipschitz}, \sup_{\mathbf{s} \in \mathcal{D}} |f(\mathbf{s})| \leq 1} \left| \int_\mathcal{D} f(\mathbf{s}) \mu(d\mathbf{s}) \right|.$$

By $f$ is 1-Lipschitz, we mean that $|f(\mathbf{s}) - f(\mathbf{s}')| \leq d(\mathbf{s},\mathbf{s}')$ for all $\mathbf{s},\mathbf{s}' \in \mathcal{D}$. It is well known that this norm induces the topology of weak convergence.

The main results are stated in the following Theorem 12. In case $\alpha < 0$, the rate of convergence depends on $I$ and it is worthwhile making the result a little more explicit. This is done, under some regular variation type hypotheses on $I$, in Corollary 13.

THEOREM 12. *The initial states $\mathbf{u}$ considered here are all deterministic.*

(i) *Suppose that $\alpha < 0$ and $\int_{l^1} \sum_{j \geq 1} s_j^{-\alpha} \mathbf{1}_{\{s_j \geq 1\}} I(d\mathbf{s}) < \infty$. Then, for every $\gamma \in [1, \Gamma]$ [$\Gamma$ is defined by (9)], there exists a positive finite constant $A$ such that, for every $\mathbf{u}$ satisfying $\sum_{j \geq 1} \exp(-u_j^\alpha) < \infty$,*

$$\|\mathcal{L}(FI^{(\mathbf{u})}(t)) - \mathcal{L}(\mathbf{U}_{\mathrm{stat}})\|$$
$$= O\left(t^{-(\gamma-1)} \int_{l^1} \sum_{j \geq 1} s_j^{-\alpha\gamma} \exp(-At^\gamma s_j^{\alpha\gamma}) I(d\mathbf{s}) + \exp(-At^\gamma u_1^{\alpha\gamma})\right)$$

*as $t \to \infty$.*

(ii) *Suppose that $\alpha = 0$ and $\int_{l^1} \sum_{j \geq 1} s_j^{1+\varepsilon} I(d\mathbf{s}) < \infty$ for some $\varepsilon > 0$. Then for every $\mathbf{u} \in l^{1+\varepsilon}$ and $a < \phi(\varepsilon)/(2+\varepsilon)$,*

$$\|\mathcal{L}(FI^{(\mathbf{u})}(t)) - \mathcal{L}(\mathbf{U}_{\mathrm{stat}})\| = o(\exp(-at)) \qquad \text{as } t \to \infty.$$

(iii) *Suppose that $\alpha > 0$ and $\int_{l^1} \sum_{j \geq 1} s_j^p I(d\mathbf{s}) < \infty$ for some $p > 0$. Then, for every $\mathbf{u} \in l^p$ and every $a < 1/\alpha$,*

$$\|\mathcal{L}(FI^{(\mathbf{u})}(t)) - \mathcal{L}(\mathbf{U}_{\mathrm{stat}})\| = o(t^{-a}) \qquad \text{as } t \to \infty.$$

Note first that, by Theorems 7, 8 and 9, the assumptions we make on $I$ imply in each case that $\mathbf{U}_{\mathrm{stat}} \in \mathcal{D}$ a.s. In case $\alpha < 0$, the given upper bound may be infinite for some $\gamma$'s. The point is then to find the $\gamma$'s in $[1, \Gamma]$ that give the best rate of convergence. This is possible, for example, when $\int_{l^1} \sum_{j \geq 1} \mathbf{1}_{\{s_j \geq x\}} I(d\mathbf{s})$ behaves regularly as $x \to \infty$. In such case the statement (i) turns to:

COROLLARY 13. *Suppose $\alpha < 0$ and fix $\mathbf{u}$ such that $\sum_{j \geq 1} \exp(-u_j^\alpha) < \infty$.*

(i) *If $\int_{l^1} \sum_{j \geq 1} \mathbf{1}_{\{s_j \geq x\}} I(d\mathbf{s}) \sim l(x) x^{-\varrho}$ as $x \to \infty$ for some slowly varying function $l$ and some $\varrho > 0$, then, provided $-\alpha < \varrho$,*

$$\|\mathcal{L}(FI^{(\mathbf{u})}(t)) - \mathcal{L}(\mathbf{U}_{\mathrm{stat}})\| = O(l(t^{1/|\alpha|}) t^{-(\varrho/|\alpha|-1)}) \qquad \text{as } t \to \infty.$$

(ii) *If $-\ln(\int_{l^1} \sum_{j \geq 1} \mathbf{1}_{\{s_j \geq x\}} I(d\mathbf{s})) \sim l(x) x^\varrho$ as $x \to \infty$ for some slowly varying function $l$ and some $\varrho > 0$, then there exists a slowly varying function $l'$ (which is constant when $l$ is constant) such that*

$$\|\mathcal{L}(FI^{(\mathbf{u})}(t)) - \mathcal{L}(\mathbf{U}_{\mathrm{stat}})\| = O(t^{-(\Gamma-1)} \exp(-l'(t) t^{\varrho \Gamma/(|\alpha|\Gamma+\varrho)})) \qquad \text{as } t \to \infty.$$

*In the special case when $I(s_1 > a) = 0$ for some $a > 0$,*

$$\|\mathcal{L}(FI^{(\mathbf{u})}(t)) - \mathcal{L}(\mathbf{U}_{\text{stat}})\| = O(\exp(-Bt^{\Gamma}))$$

*for some constant $B > 0$.*

PROOF. (i) First, by integrating by parts and then using, for example, Proposition 1.5.10 of [10], one obtains that, for $\gamma \in [1, \varrho/(-\alpha))$,

$$\int_{l^1} \sum_{j \geq 1} s_j^{-\alpha\gamma} \mathbf{1}_{\{x \geq s_j^{\alpha\gamma}\}} I(d\mathbf{s}) \approx l(x^{1/\alpha\gamma}) x^{-1-\varrho/\alpha\gamma} \qquad \text{as } x \to 0$$

(the notation $\approx$ means that the functions are equivalent up to a multiplicative constant). Then, using Karamata's Abelian–Tauberian theorem (Theorem 1.7.1′ of [10]), one deduces that

$$\int_{l^1} \sum_{j \geq 1} s_j^{-\alpha\gamma} \exp(-t s_j^{\alpha\gamma}) I(d\mathbf{s}) \approx l(t^{-1/\alpha\gamma}) t^{1+\varrho/\alpha\gamma} \qquad \text{as } t \to \infty.$$

Now if $-\alpha < \varrho$, statement (i) of Theorem 12 applies and one can plug the above equivalence into the upper bound obtained there, hence, the conclusion.

(ii) Let $1 \leq \gamma \leq \Gamma$. By integrating by parts and then by using Theorem 4.12.10 in [10], one sees that $-\ln(\int_{l^1} \sum_{j \geq 1} s_j^{-\alpha\gamma} \mathbf{1}_{\{s_j \geq x\}} I(d\mathbf{s})) \sim l(x) x^{\varrho}$ as $x \to \infty$. According to de Bruijn's Abelian–Tauberian theorem 4.12.9 in [10], this implies that

$$(19) \qquad -\ln\left(\int_{l^1} \sum_{j \geq 1} s_j^{-\alpha\gamma} \exp(-t s_j^{\alpha\gamma}) I(d\mathbf{s})\right) \approx f(t) \qquad \text{as } t \to \infty,$$

where $f(t) = 1/\Psi^{\leftarrow}(t)$ with $\Psi(t) = \Phi(t)/t$ and $\Phi^{\leftarrow}(t) = t^{\varrho/(\alpha\gamma)}/l(t^{1/(-\alpha\gamma)})$. Here $\Phi^{\leftarrow}(t) = \sup\{u \geq 0 : \phi(u) > t\}$ and similarly for $\Psi$. Therefore, $f(t) \sim \widetilde{l}(t) t^{\varrho/(\varrho+|\alpha|\gamma)}$ for some slowly varying function $\widetilde{l}$ (to invert regularly varying functions, we refer to Chapter 1.5.7 of [10]) which is constant when $l$ is constant. The assumption we have on $I$ allows us to apply Theorem 12(i) and the conclusion then follows by taking there $\gamma = \Gamma$ and using the equivalence (19). The special case when $I(s_1 > a) = 0$ is obvious. □

Hence, our bounds for the rate of convergence depend significantly on $I$ when $\alpha < 0$, whereas they are essentially independent of $I$ when $\alpha \geq 0$. Also, in any case they are essentially independent of the initial state $\mathbf{u}$.

We now turn to the proof of Theorem 12, which relies on a coupling method that holds for $\mathcal{D}$-valued $X$-processes with immigration, as defined in Section 2.1. We first explain the method in this general context and then make precise calculations for fragmentation with immigration processes. In this latter case, if $c, \nu$ and $I$ are fixed so that $I(s_1 > 1) = 0$ and if $\alpha$ varies, one sees (without any calculations, just using that particles with mass $\leq 1$ split faster when $\alpha$ is smaller) that the employed method provides a better rate of convergence when $\alpha$ is smaller. When $I(s_1 > 1) > 0$, the comparison of rates of convergence as $\alpha$ varies is no longer possible because particles with mass larger than 1 split more slowly when $\alpha$ is smaller.

PROOF OF THEOREM 12. Let $X$ be a $\mathcal{D}$-valued Markov process with the superposition property and $I$ an immigration measure such that the processes $XI^{(\mathbf{u})}$, $\mathbf{u} \in \mathcal{D}$, defined by formula (11), are $\mathcal{D}$-valued $X$-processes with immigration. Let then $((\mathbf{s}(t_i), t_i), i \geq 1)$ be the atoms of a Poisson measure with intensity $I(d\mathbf{s}) dt$, $t \geq 0$, and suppose that the stationary sum $\mathbf{U}_{\text{stat}}$ constructed from $((\mathbf{s}(t_i), t_i), i \geq 1)$, as explained in (12), belongs a.s. to $\mathcal{D}$. Suppose, moreover, that $X^{(\mathbf{u})}(t) \stackrel{\text{a.s.}}{\to} \mathbf{0}$ for all $\mathbf{u} \in \mathcal{D}$.

Then, fix $\mathbf{u} \in \mathcal{D}$ and consider $X^{(\mathbf{u})}$ and $X^{(\mathbf{U}_{\text{stat}})}$ some versions of $X$ starting, respectively, from $\mathbf{u}$ and $\mathbf{U}_{\text{stat}}$. Consider next $XI^{(\mathbf{0})}$ an $X$-process with immigration starting from $\mathbf{0}$, independent of $X^{(\mathbf{u})}$ and $X^{(\mathbf{U}_{\text{stat}})}$. Then, the processes $XI^{(\mathbf{u})}$ and $XI^{(\mathbf{U}_{\text{stat}})}$, defined, respectively, by $XI^{(\mathbf{u})}(t) := X^{(\mathbf{u})}(t) + XI^{(\mathbf{0})}(t)$ and $XI^{(\mathbf{U}_{\text{stat}})}(t) := X^{(\mathbf{U}_{\text{stat}})}(t) + XI^{(\mathbf{0})}(t)$, $t \geq 0$, are $X$-processes with immigration starting, respectively, from $\mathbf{u}$ and $\mathbf{U}_{\text{stat}}$.

Let now $r$ be a deterministic function and call $\tau_r^{(\mathbf{u})}$ the first time $t$ at which $X_1^{(\mathbf{u})}(s) \leq r(s)$ for all $s \geq t$ and, similarly, $\tau_r^{(\text{stat})}$ the first time $t$ at which $X_1^{(\mathbf{U}_{\text{stat}})}(s) \leq r(s)$ for all $s \geq t$. Of course, the interesting cases are $\tau_r^{(\mathbf{u})} < \infty$ and $\tau_r^{(\text{stat})} < \infty$ a.s. Such cases exist, take, for example, $r \equiv 1$.

Our goal is to evaluate the behavior of the norm $\|\mathcal{L}(XI^{(\mathbf{u})}(t)) - \mathcal{L}(\mathbf{U}_{\text{stat}})\|$ as $t \to \infty$. To do so, let $f : \mathcal{D} \to \mathbb{R}$ denote a 1-Lipschitz function on $\mathcal{D}$ such that $\sup_{\mathbf{s} \in \mathcal{D}} |f(\mathbf{s})| \leq 1$. For all $t \geq 0$, we construct a function $f_{r(t)}$ from $f$ and $r(t)$ by setting

$$f_{r(t)}(\mathbf{s}) := \begin{cases} f(\mathbf{0}), & \text{when } s_1 \leq r(t), \\ f(s_1, \ldots, s_{i(r(t))}, 0, 0, \ldots), & \text{when } s_1 > r(t), \end{cases}$$

where $i(r(t))$ is the unique integer such that $s_{i(r(t))} > r(t)$ and $s_{i(r(t))+1} \leq r(t)$. Clearly, as $f$ is 1-Lipschitz and $d(\mathbf{s}, \mathbf{s}') = \sup_{j \geq 1} |s_j - s'_j|$ for $\mathbf{s}, \mathbf{s}' \in \mathcal{D}$, $|f(\mathbf{s}) - f_{r(t)}(\mathbf{s})| \leq r(t)$ for every $\mathbf{s} \in \mathcal{D}$ and, therefore,

$$
\begin{aligned}
&|E[f(XI^{(\mathbf{u})}(t)) - f(\mathbf{U}_{\text{stat}})]| \\
(20) \quad &= |E[f(XI^{(\mathbf{u})}(t)) - f(XI^{(\mathbf{U}_{\text{stat}})}(t))]| \\
&\leq 2r(t) + |E[f_{r(t)}(XI^{(\mathbf{u})}(t)) - f_{r(t)}(XI^{(\mathbf{U}_{\text{stat}})}(t))]|.
\end{aligned}
$$

The time $\tau_r^{(\mathbf{u})}$ and the function $f_{r(t)}$ are defined so that, for times $t \geq \tau_r^{(\mathbf{u})}$, $f_{r(t)}(XI^{(\mathbf{u})}(t))$ takes only into account the masses of particles that are de-scended from immigrated particles, not from $\mathbf{u}$. Therefore, one has

$$E[f_{r(t)}(XI^{(\mathbf{u})}(t))] = E[f_{r(t)}(XI^{(\mathbf{u})}(t))\mathbf{1}_{\{\tau_r^{(\mathbf{u})} \vee \tau_r^{(\text{stat})} > t\}}]$$
$$+ E[f_{r(t)}(XI^{(\mathbf{0})}(t))\mathbf{1}_{\{t \geq \tau_r^{(\mathbf{u})} \vee \tau_r^{(\text{stat})}\}}]$$

and, similarly,

$$E[f_{r(t)}(XI^{(\mathbf{U}_{\text{stat}})}(t))] = E[f_{r(t)}(XI^{(\mathbf{U}_{\text{stat}})}(t))\mathbf{1}_{\{\tau_r^{(\mathbf{u})} \vee \tau_r^{(\text{stat})} > t\}}]$$
$$+ E[f_{r(t)}(XI^{(\mathbf{0})}(t))\mathbf{1}_{\{t \geq \tau_r^{(\mathbf{u})} \vee \tau_r^{(\text{stat})}\}}].$$

Combined with (20), this gives

$$|E[f(XI^{(\mathbf{u})}(t)) - f(\mathbf{U}_{\text{stat}})]|$$
$$\leq 2r(t) + |E[(f_{r(t)}(XI^{(\mathbf{u})}(t)) - f_{r(t)}(XI^{(\mathbf{U}_{\text{stat}})}(t)))\mathbf{1}_{\{\tau_r^{(\mathbf{u})} \vee \tau_r^{(\text{stat})} > t\}}]|$$
$$\leq 2r(t) + 2P(\tau_r^{(\mathbf{u})} \vee \tau_r^{(\text{stat})} > t)$$

since $\sup_{\mathbf{s} \in \mathcal{D}} |f(\mathbf{s})| \leq 1$. This holds for all 1-Lipschitz functions $f$ such that $\sup_{\mathbf{s} \in \mathcal{D}} |f(\mathbf{s})| \leq 1$ and, therefore,

$$(21) \quad \|\mathcal{L}(XI^{(\mathbf{u})}(t)) - \mathcal{L}(\mathbf{U}_{\text{stat}})\| \leq 2(r(t) + P(\tau_r^{(\mathbf{u})} > t) + P(\tau_r^{(\text{stat})} > t)).$$

The point is thus to find a function $r$ such that the above upper bound gives the best possible rate of convergence.

In the rest of this proof, we replace $X$ by an $(\alpha, c, \nu)$-fragmentation process $F$, in order to make precise calculations. We recall that $F^{(\mathbf{u})}(t) \overset{\text{a.s.}}{\to} \mathbf{0}$ and that the assumptions of Theorem 12 involving $I$ ensure that $\mathbf{U}_{\text{stat}} \in \mathcal{D}$ a.s. for all $\alpha \in \mathbb{R}$, so that (21) holds for $FI^{(\mathbf{u})}$. The choice of the function $r$ then differs according as $\alpha < 0$, $\alpha = 0$ and $\alpha > 0$.

*Proof of* (i). Here we take $r \equiv 0$. According to the definitions above, $\tau_r^{(\mathbf{u})}$ is the first time at which $F^{(\mathbf{u})}$ reaches $\mathbf{0}$ (it may be a priori infinite) and $\tau_r^{(\text{stat})}$ the first time at which $F^{(\mathbf{U}_{\text{stat}})}$ reaches $\mathbf{0}$. As recalled in Section 1.1, the first time $\tau$ at which a 1-mass particle reaches $\mathbf{0}$ is a.s. finite since $\alpha < 0$. By self-similarity, the first time at which a particle with mass $m$ is reduced to $\mathbf{0}$ is distributed as $m^{-\alpha}\tau$. Hence, by definitions of $F^{(\mathbf{u})}$ and $F^{(\mathbf{U}_{\text{stat}})}$,

$$\tau_r^{(\mathbf{u})} = \sup_{j \geq 1} u_j^{-\alpha}\tau^{(j)} \quad \text{and} \quad \tau_r^{(\text{stat})} = \sup_{i \geq 1, j \geq 1} (s_j^{-\alpha}(t_i)\tau^{(i,j)} - t_i)^+,$$

where $(\tau^{(j)}, j \geq 1)$ and $(\tau^{(i,j)}, i, j \geq 1)$ denote families of i.i.d. copies of $\tau$ such that $(\tau^{(i,j)}, i, j \geq 1)$ is independent of $((\mathbf{s}(t_i), t_i), i \geq 1)$.

Now fix $\gamma \in [1, \Gamma]$. On the one hand, one has

$$P(\tau_r^{(\mathbf{u})} > t) \leq \sum_{j \geq 1} P(\tau^{(j)} > tu_j^\alpha),$$

which by (8) is bounded from above by $C_\gamma \sum_{j\geq 1} \exp(-C'_\gamma t^\gamma u_j^{\alpha\gamma})$ for some constants $C_\gamma, C'_\gamma > 0$. Let $0 < \varepsilon < C'_\gamma$. It is easy that this sum is, in turn, bounded for all $t \geq 1$ by $B\exp(-(C'_\gamma - \varepsilon)t^\gamma u_1^{\alpha\gamma})$, where $B$ is a constant (depending on $\gamma, \varepsilon$ and $\mathbf{u}$, not on $t \geq 1$) which is finite as soon as $\sum_{j\geq 1} \exp(-u_j^\alpha) < \infty$. On the other hand,

$$P(\tau_r^{(\text{stat})} > t) \leq \int_0^\infty \int_{l^1} \sum_{j\geq 1} P(\tau > (t+v)s_j^\alpha) I(d\mathbf{s})\, dv,$$

which, again by (8), is bounded from above by

$$\frac{C_\gamma}{C'_\gamma \gamma t^{\gamma-1}} \int_{l^1} \sum_{j\geq 1} s_j^{-\alpha\gamma} \exp(-C'_\gamma t^\gamma s_j^{\alpha\gamma}) I(d\mathbf{s})$$

for $t > 0$. Hence, the result.

*Proof of* (ii). When $\alpha = 0$, the fragmentation does not reach $\mathbf{0}$ in general. We thus have to choose some function $r \neq 0$. By assumption, $\int_{l^1} \sum_{j\geq 1} s_j^{1+\varepsilon} I(d\mathbf{s}) < \infty$ for some $\varepsilon > 0$. So, fix such $\varepsilon$, fix $\eta > 1$ and set $a := \phi(\varepsilon)/(1+\eta(1+\varepsilon))$. Then take $r(t) := \exp(-at)$, $t \geq 0$.

In order to bound from above $P(\tau_r^{(\mathbf{u})} > t)$ and $P(\tau_r^{(\text{stat})} > t)$, introduce, for all $x > 0$,

$$\tau_{a,x} = \sup\{t \geq 0 : F_1(t) > x\exp(-at)\}$$

the last time $t$ at which the largest fragment of a standard fragmentation process $F$ starting from $(1, 0, \ldots)$ has a mass larger than $x\exp(-at)$. Here we use the convention $\sup(\varnothing) = 0$. This time $\tau_{a,x}$ is a.s. finite because $\exp(at)F_1(t) \overset{\text{a.s.}}{\to} 0$ when $0 \leq a < \sup_{p\geq 0} \frac{\phi(p)}{p+1}$, as explained in [8]. More precisely, one can show the existence of a positive constant $C(a)$ such that

(22) $\qquad P(\tau_{a,x} > t) \leq C(a) x^{-(1+\varepsilon)} \exp(-at) \qquad$ for all $x > 0, t \geq 1$.

Indeed, let $t \geq 1$ and note that

$$\begin{aligned}P(\eta t \geq \tau_{a,x} > t) &\leq P(\exists u \in [t, \eta t[ : F_1(u)\exp(au) > x) \\ &\leq P(F_1(t)\exp(a\eta t) > x) \qquad (\text{as } F_1 \searrow) \\ &\leq x^{-(1+\varepsilon)} \exp(a\eta(1+\varepsilon)t) E[(F_1(t))^{1+\varepsilon}].\end{aligned}$$

This last expectation is bounded from above by $E[\sum_{k\geq 1}(F_k(t))^{1+\varepsilon}] = \exp(-\phi(\varepsilon)t)$, which yields $P(\eta t \geq \tau_{a,x} > t) \leq x^{-(1+\varepsilon)} \exp(-at)$, since $a = \phi(\varepsilon) - a\eta(1+\varepsilon)$. Then, setting $C(a) := \sum_{n\geq 1} \exp(-a(\eta^{n-1} - 1))$, one obtains (22).

By definition, $\tau_r^{(\mathbf{u})}$ is the supremum of times $t$ such that $F_1^{(\mathbf{u})}(t) > \exp(-at)$. Hence, there exist some independent random variables $\tau_{a,1/u_j}^{(j)}$, $j \geq 1$, where

$\tau_{a,1/u_j}^{(j)}$ has the same distribution as $\tau_{a,1/u_j}$, such that

$$\tau_r^{(\mathbf{u})} = \sup_{j \geq 1} \tau_{a,1/u_j}^{(j)}.$$

Then, by (22),

(23) $$P(\tau_r^{(\mathbf{u})} > t) \leq C(a) \exp(-at) \sum_{j \geq 1} u_j^{1+\varepsilon}.$$

Next, by definition of $\tau_r^{(\text{stat})}$, there exists a family of r.v. $\tau_{a,\exp(at_i)/s_j(t_i)}^{(i,j)}$, $i,j \geq 1$, such that

$$\tau_r^{(\text{stat})} = \sup_{i \geq 1, j \geq 1} (\tau_{a,\exp(at_i)/s_j(t_i)}^{(i,j)} - t_i)^+$$

and, conditionally on $((\mathbf{s}(t_i), t_i), i \geq 1)$, $\tau_{a,\exp(at_i)/s_j(t_i)}^{(i,j)} \stackrel{\text{law}}{=} \tau_{a,\exp(at_i)/s_j(t_i)}$, $i, j \geq 1$, and the $\tau_{a,\exp(at_i)/s_j(t_i)}^{(i,j)}$'s are independent. This implies that

$$P(\tau_r^{(\text{stat})} > t) \leq \sum_{i \geq 1} \sum_{j \geq 1} P(\tau_{a,\exp(at_i)/s_j(t_i)}^{(i,j)} > t_i + t)$$

and then, by (22), that

$$P(\tau_r^{(\text{stat})} > t) \leq \frac{C(a)}{2a+\varepsilon} \exp(-at) \int_{l^1} \sum_{j \geq 1} s_j^{1+\varepsilon} I(d\mathbf{s}).$$

Combining this last inequality with (21) and (23), one obtains

$$\|\mathcal{L}(FI^{(\mathbf{u})}(t)) - \mathcal{L}(\mathbf{U}_{\text{stat}})\|$$
$$\leq 2\exp(-at)\left(1 + C(a) \sum_{j \geq 1} u_j^{1+\varepsilon} + (2a)^{-1} C(a) \int_{l^1} \sum_{j \geq 1} s_j^{1+\varepsilon} I(d\mathbf{s})\right).$$

This holds for every $\eta > 1$ and, therefore,

$$\|\mathcal{L}(FI^{(\mathbf{u})}(t)) - \mathcal{L}(\mathbf{U}_{\text{stat}})\| = O(\exp(-at))$$

for every $a < \phi(\varepsilon)/(2+\varepsilon)$, provided $\mathbf{u} \in l^{1+\varepsilon}$. Then, as this holds for all values of $a$ in an open interval, one can replace $O(\exp(-at))$ by $o(\exp(-at))$.

*Proof of* (iii). Fix $0 < a < 1/\alpha$ and set $r(t) := t^{-a}$, $t > 0$. By assumption, there exists some $p > 0$ such that $\int_{l^1} \sum_{j \geq 1} s_j^p I(d\mathbf{s}) < \infty$ and we call $z$ the real number such that $z\alpha^2(a+1) = p(1 - \alpha a - \alpha z)$. Note that $0 < z < \alpha^{-1} - a$. Define then, for $x > 0$,

$$\tau_{a,x} := \sup\{t \geq 0 : F_1(t) > xt^{-a}\}.$$

The fact that $z \in (0, \alpha^{-1})$ allows us to choose some $\eta > 0$ and $q > 1$ such that $\frac{q-1}{\alpha+\eta} - aq = q(\alpha^{-1} - a - z)$, which, by definition of $z$, is also equal to $qz\alpha(a+1)/p$. According to Lemma 11(ii), there exists an r.v. $D_{(\eta,q)}$ with positive moments of all orders such that

$$t^{qa} F_1^q(t) \leq D_{(\eta,q)} t^{qa - (q-1)/(\alpha+\eta)} = D_{(\eta,q)} t^{-qz\alpha(a+1)/p}$$

a.s. for every $t > 0$. This implies that

$$P(\tau_{a,x} > t) \leq P(\exists u \geq t : u^{qa} F_1^q(u) > x^q)$$
$$\leq P(\exists u \geq t : D_{(\eta,q)} u^{-qz\alpha(a+1)/p} > x^q)$$
$$\leq B x^{-p/(z\alpha)} t^{-(a+1)},$$

where $B := E[D_{(\eta,q)}^{p/(qz\alpha)}] < \infty$.

A moment of thought shows that the times $\tau_r^{(\mathbf{u})} = \sup\{t \geq 0 : F_1^{(\mathbf{u})}(t) > t^{-a}\}$ and $\tau_r^{(\text{stat})} = \sup\{t \geq 0 : F_1^{(\mathbf{U}_{\text{stat}})}(t) > t^{-a}\}$ satisfy

$$\tau_r^{(\mathbf{u})} = \sup_{j \geq 1}(u_j^{-\alpha} \tau_{a,u_j^{\alpha a - 1}}^{(j)}) \quad \text{and} \quad \tau_r^{(\text{stat})} \leq \sup_{i \geq 1, j \geq 1}(s_j^{-\alpha} \tau_{a,s_j^{\alpha a - 1}}^{(i,j)} - t_i)^+,$$

where the r.v. $\tau_{a,u_j^{\alpha a-1}}^{(j)}$, $j \geq 1$, are independent such that $\tau_{a,u_j^{\alpha a-1}}^{(j)} \stackrel{\text{law}}{=} \tau_{a,u_j^{\alpha a-1}}$ and, conditionally on $((\mathbf{s}(t_i), t_i), i \geq 1)$, the r.v. $\tau_{a,s_j^{\alpha a-1}}^{(i,j)}$, $i, j \geq 1$, are independent such that $\tau_{a,s_j^{\alpha a-1}}^{(i,j)} \stackrel{\text{law}}{=} \tau_{a,s_j^{\alpha a-1}}$. Using then the upper bound $P(\tau_{a,x} > t) \leq B x^{-p/(z\alpha)} t^{-(a+1)}$, one obtains

$$P(\tau_r^{(\mathbf{u})} > t) \leq B t^{-(a+1)} \sum_{j \geq 1} u_j^{-\alpha(a+1) + p(1-\alpha a)/z\alpha},$$

which is equal to $B t^{-(a+1)} \sum_{j \geq 1} u_j^p$ by definition of $z$. Similarly, one obtains

$$P(\tau_r^{(\text{stat})} > t) \leq a^{-1} B t^{-a} \int_{l^1} \sum_{j \geq 1} s_j^p I(d\mathbf{s}).$$

Hence, by (21),

$$\|\mathcal{L}(FI^{(\mathbf{u})}(t)) - \mathcal{L}(\mathbf{U}_{\text{stat}})\| \leq R t^{-a} \left(1 + \sum_{j \geq 1} u_j^p + \int_{l^1} \sum_{j \geq 1} s_j^p I(d\mathbf{s})\right),$$

where $R$ is a finite real number depending on the parameters of the fragmentation and on $a$, but not on $t$ and $f$. This holds for all $a \in (0, 1/\alpha)$, which gives the bounds $o(t^{-a})$, $a < 1/\alpha$, claimed in the statement. $\square$

**4. An example constructed from a Brownian motion with positive drift.**
Let $B$ be a standard linear Brownian motion and for every $d > 0$, consider the Brownian motion with drift $d$,

$$B_{(d)}(x) := B(x) + dx, \qquad x \geq 0.$$

For any $t > 0$, define

$$L_{(d)}(t) := \inf\{x \geq 0 : B_{(d)}(x) = t\} \qquad R_{(d)}(t) := \sup\{x \geq 0 : B_{(d)}(x) = t\},$$

the first and the last hitting times of $t$ by $B_{(d)}$. Clearly, $0 < L_{(d)}(t) < R_{(d)}(t) < \infty$ a.s., since $d > 0$. It is thus possible to consider the decreasing rearrangement of lengths of the connected components of

$$\mathcal{E}_{(d)}(t) := \{x \in [L_{(d)}(t), R_{(d)}(t)] : B_{(d)}(x) > t\},$$

which we denote by $FI_{(d)}(t)$.

PROPOSITION 14. (i) *The process $(FI_{(d)}(t), t \geq 0)$ is a fragmentation immigration process with the following parameters:*

- $\alpha_B = -1/2$,
- $c_B = 0$,
- $\nu_B(s_1 + s_2 < 1) = 0$ and $\nu_B(s_1 \in dx) = \sqrt{2\pi^{-1}} x^{-3/2}(1-x)^{-3/2} dx$, $x \in [1/2, 1)$,
- $I_{(d)}(s_2 > 0) = 0$ and $I_{(d)}(s_1 \in dx) = \sqrt{(2\pi)^{-1}} x^{-3/2} \exp(-xd^2/2) dx$, $x > 0$.

(ii) *The process is stationary. The stationary law is that of a Cox measure (that is, a Poisson measure with random intensity) with intensity $T(d)\sqrt{(8\pi)^{-1}} x^{-3/2} \times \exp(-xd^2/2) dx$, $x > 0$, where $T(d)$ is an exponential r.v. with parameter $d$.*

(iii) *There exists a constant $L \in (0, \infty)$ such that, for every $\mathbf{u} \in \mathcal{D}$ satisfying $\sum_{j \geq 1} \exp(-u_j^{-1/2}) < \infty$, an $(\alpha_B, c_B, \nu_B, I_{(d)})$-fragmentation immigration $FI^{(\mathbf{u})}$ starting from $\mathbf{u}$ converges in law to the stationary distribution $\mathcal{L}(\mathbf{U}_{\text{stat}})$ at rate*

$$\|\mathcal{L}(FI^{(\mathbf{u})}(t)) - \mathcal{L}(\mathbf{U}_{\text{stat}})\| = O(t^{-1} \exp(-Lt)).$$

Note that the immigrating particles arrive one-by-one.

The fragmentation part of this process, that does not depend on $d$, is a well-known fragmentation process that was first constructed by Bertoin in [7]. Let $F_B^{(\mathbf{l})}$ denote this fragmentation starting from $\mathbf{l} = (l, 0, \dots)$. It is a binary fragmentation, that is, each particle splits exactly into two pieces, which is constructed from a Brownian excursion $e_B^{(l)}$ conditioned to have length $l$ as follows:

(24) $F_B^{(\mathbf{l})}(t) := \{\text{lengths of connected components of } \{x \in [0, l] : e_B^{(l)}(x) > t\}\}^{\downarrow}$

for all $t \geq 0$. In [7] it is proved that this process is indeed a fragmentation process with index $\alpha_B = -1/2$, no erosion and a dislocation measure $\nu_B$ as given above.

PROOF OF PROPOSITION 14. (i) According to Corollaries 1 and 2 in [25], the process defined by

$$Y_{(d)}(x) := B_{(d)}(x + R_{(d)}(0)), \qquad x \geq 0,$$

is a $\text{BES}^0(3,d)$ (which means that it is identical in law to the norm of a three-dimensional Brownian motion with drift $d$) and is independent of $(B_{(d)}(x), 0 \leq x \leq R_{(d)}(0))$. This last process codes the fragmentation of particles present at time 0, whereas the process $Y_{(d)}$ codes the immigration and fragmentation of immigrated particles. More precisely:

- Let $e_B^{(l_1)}, \ldots, e_B^{(l_i)}, \ldots$ denote the finite excursions of $B_{(d)}$ above 0, with respective lengths $l_1, l_2, \ldots$. The Cameron–Martin–Girsanov theorem implies that the $(l_i, i \geq 1)$ are the finite jumps of a subordinator with Lévy measure $\sqrt{(8\pi)^{-1}} x^{-3/2} e^{-xd^2/2} dx$, killed at an exponential time with parameter $d$, and that conditionally on $(l_i, i \geq 1)$ the excursions $e_B^{(l_1)}, e_B^{(l_2)}, \ldots$ are independent Brownian excursions with respective lengths $l_1, \ldots, l_i, \ldots$. This gives the distribution of $FI_{(d)}(0) = (l_1, l_2, \ldots)^{\downarrow}$ and implies that the process $(FI_{(d)}^{[0,R_{(d)}(0)]}(t), t \geq 0)$ defined by

$$FI_{(d)}^{[0,R_{(d)}(0)]}(t) := \{\text{lengths of connected}$$
$$\text{comp. of } \{x \in [L_{(d)}(t), R_{(d)}(0)] : B_{(d)}(x) > t\}\}^{\downarrow}$$

is an $(-1/2, 0, \nu_B)$-fragmentation starting from $FI_{(d)}(0)$.

- Let $J_{(Y_{(d)})}(x) := \inf_{y \geq x} Y_{(d)}(y)$, $x \geq 0$, be the future infimum of $Y_{(d)}$. One has to see $J_{(Y_{(d)})}$ as the process coding the arrival of immigrating particles and $Y_{(d)} - J_{(Y_{(d)})}$ as the process coding their fragmentation. According to a generalization of Pitman's theorem (Corollary 1 in [25]), $(J_{(Y_{(d)})}, Y_{(d)} - J_{(Y_{(d)})})$ is distributed as $(M_{(d)}, M_{(d)} - B_{(d)})$, where $M_{(d)}(x) := \sup_{[0,x]} B_{(d)}(y)$, $x \geq 0$. Moreover, according to the Cameron–Martin–Girsanov theorem, $M_{(d)}$ is distributed as the inverse of a subordinator with Lévy measure

$$I_{(d)}(s_1 \in dx) = \sqrt{(2\pi)^{-1}} x^{-3/2} \exp(-xd^2/2) \, dx, \qquad x > 0,$$

and, conditionally on their lengths, the excursions above 0 of $M_{(d)} - B_{(d)}$ are Brownian excursions. Let $((\Delta_{(d)}(t_i), t_i), i \geq 1)$ denote the family of jump sizes and times of the subordinator inverse of $M_{(d)}$. The sequence

$$FI_{(d)}^{[R_{(d)}(0),\infty)}(t) := \{\text{lengths of connected}$$
$$\text{comp. of } \{x \in [R_{(d)}(0), R_{(d)}(t)] : B_{(d)}(x) > t\}\}^{\downarrow}$$

is the decreasing rearrangement of masses of particles that have immigrated at time $t_i \leq t$ with mass $\Delta_{(d)}(t_i)$ and that have split independently (conditionally on their masses) until time $t - t_i$ according to the fragmentation $(-1/2, 0, \nu_B)$.

- $FI_{(d)}(t)$ is the concatenation of $FI_{(d)}^{[0,R_{(d)}(0)]}(t)$ and $FI_{(d)}^{[R_{(d)}(0),\infty)}(t)$, which leads to the result. Note that $I_{(d)}$ satisfies the hypothesis (H1).

(ii) That $FI_{(d)}(t) \stackrel{\text{law}}{=} FI_{(d)}(0)$ is a simple consequence of the strong Markov property of $B$ applied at time $L_{(d)}(t)$. The stationary distribution $\mathcal{L}(FI_{(d)}(0))$ is calculated in the first part of this proof.

(iii) It is easy to check that the $\nu_B$-dependent parameter $\Gamma_B$ [defined in (9)] is here equal to 2 and that

$$-\ln\left(\int_x^\infty I_{(d)}(s_1 \in dy)\right) \sim \frac{d^2 x}{2} \qquad \text{as } x \to \infty.$$

Then we conclude with Corollary 13(ii). □

REMARK. Let $Y_{(d)}$ be a $\text{BES}^0(3, d)$, $d \geq 0$, and set

$FI_{Y_{(d)}}(t) := \{\text{lengths of connected}$

$\quad\text{comp. of } \{x \in [L_{Y_{(d)}}(t), R_{Y_{(d)}}(t)] : Y_{(d)}(x) > t\}\}^\downarrow,$

where $L_{Y_{(d)}}(t) := \inf\{x \geq 0 : Y_{(d)}(x) = t\}$ and $R_{Y_{(d)}}(t) := \sup\{x \geq 0 : Y_{(d)}(x) = t\}$. According to the proof above, $FI_{Y_{(d)}}$ is an $(-1/2, 0, \nu_B, I_{(d)})$-fragmentation with immigration starting from $\mathbf{0}$ (clearly, this is also valid for $d = 0$). Recall then the construction of the stationary state $\mathbf{U}_{\text{stat}}$, as explained in (12). It is easy to see that $\mathbf{U}_{\text{stat}}$ has the same law as the point measure whose atoms are the lengths of the excursions below 0 of the process obtained by reflecting $Y_{(d)}$ at the level of its future infimum. By Corollary 1 in [25], this reflected process is a Brownian motion with drift $d$. Therefore, if $d > 0$, $\mathbf{U}_{\text{stat}} \in \mathcal{D}$ a.s. and the stationary distribution is that of the reordering of the lengths of the excursions below 0 of a Brownian motion with drift $d$, which is indeed the distribution of $FI_{(d)}(0)$ (by Girsanov's theorem). On the other hand, if $d = 0$, $\mathbf{U}_{\text{stat}}$ is clearly not in $\mathcal{D}$ a.s. and then there is no stationary distribution [which was already known, according to Theorem 7(ii)].

At last, we mention that one can construct in a similar way some fragmentation with immigration processes from height functions coding continuous state branching processes with immigration (as introduced by Lambert [22]). This is detailed in [20].

**5. The fragmentation with immigration equation.** The deterministic counterpart of the fragmentation with immigration process $(\alpha, c, \nu, I)$ is the following equation, namely, the *fragmentation with immigration equation* $(\alpha, c, \nu, I)$

$$
\text{(E)} \quad \begin{aligned} \partial_t \langle \mu_t, f \rangle &= \int_0^\infty x^\alpha \bigg( -cxf'(x) + \int_{\mathcal{D}_1} \bigg[ \sum_{j \geq 1} f(xs_j) - f(x) \bigg] \nu(d\mathbf{s}) \bigg) \mu_t(dx) \\ &\quad + \int_{l^1} \sum_{j \geq 1} f(s_j) I(d\mathbf{s}), \end{aligned}
$$

where $(\mu_t, t \geq 0)$ is a family of nonnegative Radon measures on $(0, \infty)$. The measure $\mu_t(dx)$ corresponds to the average number per unit volume of particles with mass in the interval $(x, x + dx)$ at time $t$. The test-functions $f$ belong to $\mathcal{C}_c^1(0, \infty)$, the set of continuously differentiable functions with compact support in $(0, \infty)$. Note that the hypothesis (H1) implies the finiteness of the integral $\int_{l^1} \sum_{j \geq 1} f(s_j) I(d\mathbf{s})$ for every $f \in \mathcal{C}_c^1(0, \infty)$. In [2], the stationary solution to this equation is studied in the special case when $\alpha = 1$, $c = 0$, $\nu(s_1 \in dx) = 2\mathbf{1}_{\{x \in [1/2, 1]\}} dx$ and $\nu(s_1 + s_2 < 1) = 0$, $I(s_2 > 0) = 0$ and $I(s_1 \in dx) = i(x) dx$ for some measurable function $i$. Here we investigate solutions and stationary solutions to (E) in the general case.

5.1. *Solutions to* (E). When $I = 0$, existence and uniqueness of a solution to (E) starting from $\delta_1(dx)$ are established in Theorem 3 in [18]. More precisely, the unique solution to the equation starting from $\delta_1(dx)$ is given, for all $t \geq 0$, by

$$
(25) \qquad \langle \eta_t, f \rangle := E \bigg[ \sum_{k \geq 1} f(F_k(t)) \bigg], \qquad f \in \mathcal{C}_c^1(0, \infty),
$$

where $F$ is a standard fragmentation process $(\alpha, c, \nu)$. Now, we generalize this to the case when $I \neq 0$. In that aim, we recall that some fragmentation with immigration processes starting from $\mathbf{u} \in \mathcal{R}$ were introduced in (10). Recall also that $\phi$ is the Laplace exponent given by (3) and that $\overline{\phi} = \phi - \phi(0)$.

PROPOSITION 15. *Let $\mu_0$ be a nonnegative Radon measure on $(0, \infty)$ and let $\mathbf{u}$ be a Poisson measure with intensity $\mu_0$. Consider then an $(\alpha, c, \nu, I)$-fragmentation with immigration $(FI^{(\mathbf{u})}(t), t \geq 0)$, as introduced in (10), and define a family of nonnegative measures $(\mu_t, t \geq 0)$ by*

$$
(26) \qquad \langle \mu_t, f \rangle := E \bigg[ \sum_{k \geq 1} f(FI_k^{(\mathbf{u})}(t)) \bigg], \qquad f \in \mathcal{C}_c^1(0, \infty), \ f \geq 0.
$$

*If one of the three following assertions is satisfied:*

(A1) $\alpha > 0$, $\int_{l^1} \sum_{j\geq 1} s_j I(d\mathbf{s}) < \infty$ and $\int_1^\infty x \mu_0(dx) < \infty$,

(A2) $\alpha = 0$, $\int_{l^1} \sum_{j\geq 1} s_j \overline{\phi}(\frac{1}{\ln s_j}) \mathbf{1}_{\{s_j \geq 1\}} I(d\mathbf{s}) < \infty$ and $\int_1^\infty x \overline{\phi}(\frac{1}{\ln x}) \mu_0(dx) < \infty$,

(A3) $\alpha < 0$, $\int_{l^1} \sum_{j\geq 1} s_j^{1+\alpha} \mathbf{1}_{\{s_j \geq 1\}} I(d\mathbf{s}) < \infty$ and $\int_1^\infty x^{1+\alpha} \mu_0(dx) < \infty$,

then the measures $\mu_t$, $t \geq 0$, are Radon and the family $(\mu_t, t \geq 0)$ is the unique solution to the fragmentation with immigration equation (E) starting from $\mu_0$.

Of course, $FI^{(\mathbf{u})}$ is a "usual" $\mathcal{D}$-valued fragmentation with immigration process as soon as $\mu_0[1, \infty) < \infty$.

REMARKS. 1. Notice that, for all $f \in \mathcal{C}_c^1(0, \infty)$, $f \geq 0$,

$$\langle \mu_t, f \rangle = E\bigg[\sum_{i\geq 1} \sum_{k\geq 1} f(u_i F_k(u_i^\alpha t))\bigg]$$
$$+ E\bigg[\sum_{t_i \leq t} \sum_{j\geq 1} \sum_{k\geq 1} f(s_j(t_i) F_k(s_j^\alpha(t_i)(t - t_i)))\bigg],$$

where $((\mathbf{s}(t_i), t_i), i \geq 1)$ [resp. $(u_i, i \geq 1)$] are the atoms of a Poisson measure with intensity $I(d\mathbf{s})\, dt$ (resp. $\mu_0$) and $F$ is an $(\alpha, c, \nu)$-fragmentation, independent of these Poisson measures. By (6), this can be written as

$$\langle \mu_t, f \rangle = \int_0^\infty E[f(x \exp(-\xi(\rho(x^\alpha t)))) \exp(\xi(\rho(x^\alpha t)))] \mu_0(dx)$$
(27)
$$+ \int_0^t \int_{l^1} \sum_{j\geq 1} E[f(s_j \exp(-\xi(\rho(s_j^\alpha u)))) \exp(\xi(\rho(s_j^\alpha u)))] I(d\mathbf{s})\, du,$$

where $\xi$ is a subordinator with Laplace exponent $\phi$. It is not hard to see that there exist some dislocation measures $\nu_1 \neq \nu_2$ that lead to the same $\phi$. In this case, the previous formula shows that the $(\alpha, c, \nu_1, I)$- and $(\alpha, c, \nu_2, I)$-fragmentation with immigration equations have identical solutions.

2. Assume that one of the assertions (A1), (A2) and (A3) is satisfied, so that the measures $\mu_t, t \geq 0$, are Radon. Then, these measures are hydrodynamic limits of fragmentation with immigration processes. Indeed, let $\mathbf{u}^{(n)}$ be a Poisson measure with intensity $n\mu_0$ and call $FI^{(n)}$ a fragmentation with immigration process with parameters $(\alpha, c, \nu, nI)$ starting from $\mathbf{u}^{(n)}$. Then, for every $t \geq 0$,

$$\frac{1}{n} FI^{(n)}(t) \stackrel{\text{vaguely}}{\to} \mu_t(dx) \qquad \text{a.s.}$$

This holds because $FI^{(n)}(t)$ is the sum of $n$ i.i.d. point measures distributed as $FI^{(\mathbf{u}^{(1)})}(t)$ for some $(\alpha, c, \nu, I)$-fragmentation with immigration $FI^{(\mathbf{u}^{(1)})}$.

The strong law of large numbers then implies that, for every $f \in \mathcal{C}_c^1(0, \infty)$,

$$\frac{1}{n} \sum_{k \geq 1} f(FI_k^{(n)}(t)) \stackrel{\text{a.s.}}{\to} E\left[\sum_{k \geq 1} f(FI_k^{(\mathbf{u}^{(1)})}(t))\right] = \langle \mu_t, f \rangle$$

and the conclusion follows by inverting the order of "for every $f \in \mathcal{C}_c^1(0, \infty)$" and "a.s.," which can be done, for example, as in the proof of Corollary 5 of [18].

PROOF OF PROPOSITION 15. Let $\mu_t, t \geq 0$, be defined by (27) [equivalently, (26)].

• It is easily seen that these measures are Radon if (A1) holds. To prove this is also valid for assertions (A2) or (A3), we need to evaluate the rate of convergence to 0 of $P(a \leq x \exp(-\xi(\rho(x^\alpha t))) \leq b)$ as $x \to \infty$, $0 < a < b < \infty$, when $\alpha \leq 0$. First, note that this probability is bounded from above by $P(x \exp(-\overline{\xi}(\rho(x^\alpha t))) \leq b)$, where $\overline{\xi} = \xi \mathbf{1}_{\{\xi < \infty\}}$ is a subordinator with Laplace exponent $\overline{\phi} = \phi - \phi(0)$. Then for $u \geq 0$ and $v > 0$,

$$(28) \quad P(\overline{\xi}(u) > v) \leq (1 - e^{-1})^{-1} E[1 - \exp(-v^{-1}\overline{\xi}(u))]$$
$$= (1 - e^{-1})^{-1}(1 - \exp(-u\overline{\phi}(v^{-1}))).$$

When $\alpha = 0$, this implies that

$$(29) \quad P(a \leq x \exp(-\xi(t)) \leq b) = O(\overline{\phi}((\ln x)^{-1})) \quad \text{as } x \to \infty.$$

When $\alpha < 0$, by the definition of $\rho$ and conditionally on $2x^\alpha t \leq \rho(x^\alpha t) < \infty$,

$$2x^\alpha t \exp(\alpha \overline{\xi}(2x^\alpha t)) \leq \int_0^{2x^\alpha t} \exp(\alpha \overline{\xi}(r)) \, dr \leq \int_0^{\rho(x^\alpha t)} \exp(\alpha \overline{\xi}(r)) \, dr = x^\alpha t$$

and, consequently, $P(2x^\alpha t \leq \rho(x^\alpha t) < \infty) \leq P(\exp(\alpha \overline{\xi}(2x^\alpha t)) \leq 1/2)$ which, by (28), is an $O(x^\alpha)$ as $x \to \infty$. Moreover, again by (28), $P(x \exp(-\overline{\xi}(2x^\alpha t)) \leq b) = O(x^\alpha)$ and, therefore,

$$(30) \quad P(a \leq x \exp(-\xi(\rho(x^\alpha t))) \leq b) = O(x^\alpha) \quad \text{as } x \to \infty$$

since

$$P(a \leq x \exp(-\xi(\rho(x^\alpha t))) \leq b)$$
$$\leq P(2x^\alpha t \leq \rho(x^\alpha t) < \infty) + P(x \exp(-\overline{\xi}(2x^\alpha t)) \leq b).$$

Now, suppose that (A2) or (A3) holds and take $f(x) = x\mathbf{1}_{\{x \in (a,b)\}}$, $0 < a < b < \infty$. Using the results (29) and (30), one sees that $\langle \mu_t, f \rangle$ is finite. Hence, $\mu_t$ is Radon.

• Suppose that (A1), (A2) or (A3) holds, so that the measures $\mu_t, t \geq 0$, are Radon. Consider then the measures $\eta_t$, $t \geq 0$, introduced in (25). One checks that

$$\langle \mu_t, f \rangle = \int_0^\infty \langle \eta_{x^\alpha t}, f_x \rangle \mu_0(dx) + \int_0^t \int_{l^1} \sum_{j \geq 1} \langle \eta_{s_j^\alpha u}, f_{s_j} \rangle I(d\mathbf{s}) \, du,$$

where $f_x : y \mapsto f(xy)$, $x \in (0, \infty)$, $f \in \mathcal{C}_c^1(0, \infty)$. Theorem 3 in [18] states that $(\eta_t, t \geq 0)$ is a solution to (E) when $I = 0$, that is,

$$\langle \eta_t, f \rangle = f(1) + \int_0^t \langle \eta_v, Af \rangle \, dv,$$

where

$$(31) \qquad Af(x) = x^\alpha \left( -cxf'(x) + \int_{\mathcal{D}_1} \left[ \sum_{j \geq 1} f(xs_j) - f(x) \right] \nu(d\mathbf{s}) \right).$$

This equation relies on the fact that, for $f \in \mathcal{C}_c^1(0, \infty)$, $A(\mathrm{id} \times f)(x) = x^{1+\alpha} G(f)(x)$, where $G$ is the infinitesimal generator of the process $\exp(-\xi)$ (see the proof of Theorem 3 in [18] for details).

Using then that $x^\alpha A f_x = (Af)_x$, one obtains

$$(32) \qquad \langle \eta_{x^\alpha t}, f_x \rangle = f(x) + \int_0^t \langle \eta_{x^\alpha v}, (Af)_x \rangle \, dv$$

and, therefore, by Fubini's theorem,

$$\langle \mu_t, f \rangle = \langle \mu_0, f \rangle + \int_0^t \int_0^\infty \langle \eta_{x^\alpha u}, (Af)_x \rangle \mu_0(dx) \, du$$

$$+ \int_0^t \left( \int_0^u \int_{l^1} \sum_{j \geq 1} \langle \eta_{s_j^\alpha v}, (Af)_{s_j} \rangle I(d\mathbf{s}) \, dv + \int_{l^1} \sum_{j \geq 1} f(s_j) I(d\mathbf{s}) \right)$$

$$= \langle \mu_0, f \rangle + \int_0^t \langle \mu_u, Af \rangle \, du + t \int_{l^1} \sum_{j \geq 1} f(s_j) I(d\mathbf{s}).$$

(to see why Fubini's theorem holds, call $[a, b]$ the support of $f$ and suppose $f \geq 0$. The same argument holds for the integral involving $I$). Hence, $(\mu_t, t \geq 0)$ is indeed a solution to (E). It remains to prove the uniqueness. This can be done with some minor changes by adapting the proof of uniqueness of a solution to (E) when $I = 0$ (see the third part of the proof of Theorem 3 in [18]). □

5.2. *Stationary solutions to* (E). As in the stochastic case, we are interested in the existence of a stationary regime. We say that a Radon measure $\mu_{\text{stat}}$ is a stationary solution to (E) if the family $(\mu_t = \mu_{\text{stat}}, t \geq 0)$ is a solution to (E).

PROPOSITION 16. (i) *There is a stationary solution to* (E) *as soon as* $\int_{l^1} \sum_{j\geq 1} s_j I(d\mathbf{s}) < \infty$ *and, conversely, provided that hypothesis* (H2) *holds, there is no stationary solution to* (E) *when* $\int_{l^1} \sum_{j\geq 1} s_j I(d\mathbf{s}) = \infty$. *In case* $\int_{l^1} \sum_{j\geq 1} s_j I(d\mathbf{s}) < \infty$, *the stationary solution* $\mu_{\text{stat}}$ *is unique and given by*

$$\mu_{\text{stat}}(dx) := x^{-\alpha} \mu_{\text{stat}}^{(\text{hom})}(dx), \qquad x \geq 0,$$

*where the measure* $\mu_{\text{stat}}^{(\text{hom})}$ *is independent of* $\alpha$ *and is constructed from* $c, \nu$ *and* $I$ *by*

$$(33) \quad \langle \mu_{\text{stat}}^{(\text{hom})}, f \rangle := \int_0^\infty \int_{l^1} \sum_{j\geq 1} E[f(s_j \exp(-\xi(t))) \exp(\xi(t))] I(d\mathbf{s}) \, dt,$$

$$f \in \mathcal{C}_c^1(0, \infty).$$

(ii) *Suppose* $\int_{l^1} \sum_{j\geq 1} s_j I(d\mathbf{s}) < \infty$ *and* $\int_1^\infty x \mu_0(dx) < \infty$ *and let* $(\mu_t, t \geq 0)$ *be the solution to* (E) *starting from* $\mu_0$. *Then,*

$$\mu_t \stackrel{\text{vaguely}}{\to} \mu_{\text{stat}} \qquad as \ t \to \infty.$$

REMARKS. 1. It $\mu_{\text{stat}}$ exists, then $\mathbf{U}_{\text{stat}} \in \mathcal{R}$ a.s. and $\langle \mu_{\text{stat}}, f \rangle = E[\langle \mathbf{U}_{\text{stat}}, f \rangle]$, $f \in \mathcal{C}_c^1(0, \infty)$. Note that it is possible that $\mathbf{U}_{\text{stat}} \in \mathcal{R} \setminus \mathcal{D}$, which then implies that there exists no stationary solution in the stochastic case, although there is one in the deterministic case. Conversely, $\mathbf{U}_{\text{stat}}$ may belong to $\mathcal{D}$ a.s., even if its "expectation" measure $\mu$ defined by $\langle \mu, f \rangle := E[\langle \mathbf{U}_{\text{stat}}, f \rangle]$ is not Radon. Then there exists a stationary solution in the stochastic case, but not in the deterministic one.

2. Call $\Lambda := \sup\{\lambda : \int_{l^1} \sum_{j\geq 1} s_j^\lambda I(d\mathbf{s}) < \infty\}$ and suppose $\Lambda > 1$. Then the statement (i) and the relations $E[e^{-q\xi(t)}] = e^{-t\phi(q)}$, $t, q \geq 0$, imply that, for all $1 + \alpha < \lambda < \Lambda + \alpha$,

$$(34) \quad \int_0^\infty x^\lambda \mu_{\text{stat}}(dx) = \phi(\lambda - \alpha - 1)^{-1} \int_{l^1} \sum_{j\geq 1} s_j^{\lambda-\alpha} I(d\mathbf{s}),$$

and that this integral is infinite as soon as $\lambda > \Lambda + \alpha$ or $\lambda \leq 1 + \alpha$, provided $\phi(0) = 0$ [which is equivalent to $c = \nu(\sum_{j\geq 1} s_j < 1) = 0$]. This characterizes $\mu_{\text{stat}}$ and is more explicit than (33).

As an example, it allows us to obtain the more convenient expression

$$\mu_{\text{stat}}(dx) = \left( x^{-\alpha} i(x) + 2 x^{-\alpha-2} \int_x^\infty y i(y) \, dy \right) dx$$

in case $\nu$ is binary, $\nu(s_1 \in dx) = 2 \mathbf{1}_{\{x \in [1/2, 1]\}} \, dx$, $c = 0$, and $I(s_1 \in dx) = i(x) \, dx$, $I(s_2 > 0) = 0$ ($\alpha \in \mathbb{R}$). This latter result is proved in a different way in [2].

Others examples are given by the equations corresponding to the fragmentation with immigration processes constructed from Brownian motions with drift $d > 0$ (Section 4). The immigration measure $I_{(d)}$ satisfies $\int_{l^1} \sum_{j \geq 1} s_j^\lambda I_{(d)}(d\mathbf{s}) < \infty$ for all $\lambda > 1/2$ and, therefore, there exists a stationary solution to the equation. One can use (34) to obtain

$$\mu_{\text{stat}}(dx) = \frac{1}{d\sqrt{8\pi x^3}} \exp(-xd^2/2)\, dx, \qquad x \geq 0.$$

This can also be shown by using remark 1 above and the stationary law $\mathcal{L}(\mathbf{U}_{\text{stat}})$ given in Proposition 14(ii).

PROOF OF PROPOSITION 16. (i) We first suppose that there exists a stationary solution $\mu_t = \mu_{\text{stat}}$, $t \geq 0$, to (E). Of course, then $\partial_t \langle \mu_t, f \rangle = 0$ for every $t \geq 0$ and $f \in \mathcal{C}_c^1(0, \infty)$, and, consequently,

$$\langle \mu_{\text{stat}}, Af \rangle = -\int_{l^1} \sum_{j \geq 1} f(s_j) I(d\mathbf{s}),$$

where $Af$ is given by (31). Letting $t \to \infty$ in (32), we get by dominated convergence that $\langle \eta_{x^\alpha t}, f_x \rangle \to 0$ and then that $f(x) = -\int_0^\infty \langle \eta_{x^\alpha v}, (Af)_x \rangle\, dv$, $x \in (0, \infty)$. Hence,

$$\langle \mu_{\text{stat}}, Af \rangle = \int_{l^1} \sum_{j \geq 1} \int_0^\infty \langle \eta_{s_j^\alpha v}, (Af)_{s_j} \rangle\, dv\, I(d\mathbf{s}).$$

We point out that this formula characterizes $\mu_{\text{stat}}$, since $A(\text{id} \times f)(x) = x^{1+\alpha} G(f)(x)$, where $G$ is the infinitesimal generator of $\exp(-\xi)$ and since $G(\mathcal{C}_c^1(0, \infty))$ is dense in the set of continuous functions on $(0, \infty)$ that vanish at 0 and $\infty$. Using then the definition of $\eta_t$ and formula (6), one sees that, for every measurable function $g$ with compact support in $(0, \infty)$,

$$\langle \mu_{\text{stat}}, g \rangle = \int_{l^1} \sum_{j \geq 1} \int_0^\infty E[g(s_j \exp(-\xi(\rho(s_j^\alpha v)))) \exp(\xi(\rho(s_j^\alpha v)))]\, dv\, I(d\mathbf{s})$$

(35)

$$= \int_{l^1} \sum_{j \geq 1} s_j^{-\alpha} \int_0^\infty E[g(s_j \exp(-\xi(v))) \exp((1+\alpha)\xi(v))]\, dv\, I(d\mathbf{s}),$$

using for the last equality the change of variables $v \mapsto \rho(s_j^\alpha v)$ and that $\exp(\alpha \xi_{\rho(v)})\, d\rho(v) = dv$ on $[0, D)$, $D = \inf\{v : \xi_{\rho(v)} = \infty\}$. This gives the required expression for $\mu_{\text{stat}}$.

Note now that the previous argument implies that a stationary solution exists if and only if

$$\int_{l^1} \sum_{j \geq 1} \int_0^\infty E[g(s_j \exp(-\xi(v))) \exp(\xi(v))]\, dv\, I(d\mathbf{s}) < \infty$$

for all functions $g$ of type $g(x) = x\mathbf{1}_{\{a \leq x \leq b\}}$, $0 < a < b$. For such function $g$, the previous integral is equal to

$$(36) \qquad \int_{l^1} \sum_{j \geq 1} s_j \mathbf{1}_{\{s_j \geq a\}} E[T^\xi_{\ln(s_j/a)} - T^\xi_{\ln^+(s_j/b)}] I(d\mathbf{s}),$$

where $T^\xi_t := \inf\{u : \xi(u) > t\}$, $t \geq 0$. If hypothesis (H2) holds and $\xi$ is arithmetic [i.e., if (H3) holds], the renewal theorem applies (see, e.g., Theorem I.21 in [5]) and $E[T^\xi_{\ln(t/a)} - T^\xi_{\ln^+(t/b)}]$ converges as $t \to \infty$ to some finite nonzero limit. In such case, the integral (36) is finite if and only if $\int_{l^1} \sum_{j \geq 1} s_j \mathbf{1}_{\{s_j \geq 1\}} I(d\mathbf{s}) < \infty$, $\forall b > a > 0$, and, therefore, there exists a stationary solution if and only if $\int_{l^1} \sum_{j \geq 1} s_j \mathbf{1}_{\{s_j \geq 1\}} I(d\mathbf{s}) < \infty$. This conclusion remains valid if (H2) holds and $\xi$ is not arithmetic, since the renewal theory then implies that $\limsup_{t \to \infty} E[T^\xi_{\ln(t/a)} - T^\xi_{\ln^+(t/b)}] < \infty$, and that $\liminf_{t \to \infty} E[T^\xi_{\ln(t/a)} - T^\xi_{\ln^+(t/b)}] > 0$ as soon as $\ln b - \ln a$ is large enough.

Last, to conclude when (H2) does not hold, remark first that $T^\xi_t = T^{\overline{\xi}}_t \wedge \mathbf{e}(k)$ [the subordinator $\overline{\xi}$ and the exponential r.v. $\mathbf{e}(k)$ are those defined in Section 1.1] and then that

$$E[T^\xi_{\ln(s_j/a)} - T^\xi_{\ln^+(s_j/b)}] \leq E[T^{\overline{\xi}}_{\ln(s_j/a)} - T^{\overline{\xi}}_{\ln^+(s_j/b)}] \leq E[T^{\overline{\xi}}_{\ln(b/a)}] < \infty.$$

In this case, the integral (36) is finite as soon as $\int_{l^1} \sum_{j \geq 1} s_j \mathbf{1}_{\{s_j \geq 1\}} I(d\mathbf{s}) < \infty$, $\forall b > a > 0$.

(ii) Under the assumptions of the statement, the measures $\mu_t$, $t \geq 0$, are Radon and therefore satisfy (27) for all continuous function $f$ with compact support in $(0, \infty)$. The integral involving $\mu_0$ converges to 0 as $t \to \infty$, since, with the assumption $\int_1^\infty x \mu_0(dx) < \infty$, the dominated convergence theorem applies. Hence, $\langle \mu_t, f \rangle \underset{t \to \infty}{\to} \langle \mu_{\text{stat}}, f \rangle$, using the definition (35) of $\mu_{\text{stat}}$. $\square$

**Acknowledgment.** I would like to thank my advisor, Jean Bertoin, for his constant help and support.

Laboratoire de Probabilités
et Modéles Aléatoires
Université Pierre et Marie Curie
et C.N.R.S. UMR 7599
175 rue du Chevaleret
F-75013 Paris
France
e-mail: haas@ccr.jussieu.fr